\pdfoutput=1 

\newif\REVIEW

\ifdefined\REVIEW
    \documentclass[
        a4paper,
        onecolumn
    ]{article}
    \usepackage{lineno, setspace}
\else
    \documentclass[
        a4paper,
        twocolumn
    ]{article}
\fi

\usepackage{IEK10} 
\usepackage{natbib}

\usepackage{upgreek}
\usepackage{bbm}

\DeclareSIUnit\year{yr}
\DeclareSIUnit\annum{a}
\DeclareSIUnit\ton{t}
\DeclareSIUnit{\million}{\text{M}}

\def\IEK10{
  Institute of Climate and Energy Systems,
  Energy Systems Engineering (ICE-1),
  Forschungszentrum J\"ulich GmbH,
  J\"ulich 52425,
  Germany
}
\def\RWTH{
  RWTH Aachen University,
  Aachen 52062,
  Germany
}
\def\JARA{
  JARA-ENERGY,
  J{\"u}lich 52425,
  Germany
}
\def\SVT{
  RWTH Aachen University,
  Process Systems Engineering (AVT.SVT),
  Aachen 52074,
  Germany
}

\newcommand{\mytitle}{Robust Energy System Design via Semi-infinite Programming}

\newcommand{\affil}{
  \begin{itemize}[leftmargin=3mm, itemsep=0mm]
    \item[$^a$]\IEK10
    \item[$^b$]\RWTH
    \item[$^c$]\SVT
    \item[$^d$]\JARA
  \end{itemize}
}

\def\firstAuthor{Moritz Wedemeyer}
\newcommand{\myauthor}{\firstAuthor$^{a,b}$,
Eike Cramer$^c$,
Alexander Mitsos$^{d,a,c}$,
Manuel Dahmen$^{a,*}$
}
\newcommand{\correspondingauthor}{M. Dahmen, \IEK10 \newline E-mail: m.dahmen@fz-juelich.de}
\author{\myauthor}

\usepackage[
  colorlinks,
  linkcolor=blue,
  citecolor=blue,
  urlcolor=blue,
  pdftitle={\mytitle},
  pdfauthor={\firstAuthor}
]{hyperref}
\usepackage[capitalise, nameinlink]{cleveref}
\crefname{table}{Tab.}{Tab.}

\usepackage[colorinlistoftodos]{todonotes}

\begin{document}
\twocolumn[
\begin{@twocolumnfalse}

  \thispagestyle{firststyle}

  \begin{center}
    \begin{large}
      \textbf{\mytitle}
    \end{large} \\
    \myauthor
  \end{center}

  \vspace{0.5cm}

  \begin{footnotesize}
    \affil
  \end{footnotesize}

  \vspace{0.5cm}
  
    Time-series information needs to be incorporated into energy system optimization to account for the uncertainty of renewable energy sources.
    Typically, time-series aggregation methods are used to reduce historical data to a few representative scenarios but they may neglect extreme scenarios, which disproportionally drive the costs in energy system design.
    We propose the robust energy system design (RESD) approach based on semi-infinite programming and use an adaptive discretization-based algorithm to identify worst-case scenarios during optimization.
    The RESD approach can guarantee robust designs for problems with nonconvex operational behavior, which current methods cannot achieve.
    The RESD approach is demonstrated by designing an energy supply system for the island of La Palma.
    To improve computational performance, principal component analysis is used to reduce the dimensionality of the uncertainty space.
    The robustness and costs of the approximated problem with significantly reduced dimensionality approximate the full-dimensional solution closely.
    Even with strong dimensionality reduction, the RESD approach is computationally intense and thus limited to small problems.

  \vspace{0.5cm}

  \noindent \textbf{Keywords}: \textit{Semi-infinite programming, Robust optimization, Energy system design, Variable renewable energy sources}

  \vspace{0.5cm}

  \vspace*{5mm}

\end{@twocolumnfalse}
]
\ifdefined\REVIEW
  \onecolumn
\else
  \newpage
\fi

\section{Introduction}\label{Intro}
The transformation of energy supply systems from fossil energy sources to variable renewable energy resources (VRES) has increased the volatility of electricity supply and, consequently, demand.
This volatility introduced by the stochastic nature of VRES has led to new challenges in the design and operation of energy systems \citep{grossRenewablesGridUnderstanding2007, huberIntegrationWindSolar2014, heptonstallSystematicReviewCosts2021}.

Mathematical optimization is an effective tool to design energy systems that are optimal, e.g., have minimal total annualized cost or global warming impact, and can be leveraged to design energy systems that are robust towards the volatility introduced by VRES \citep{bieglerRetrospectiveOptimization2004, yuntDesigningManportablePower2008,  zavala, liOptimalDesignOperation2015}.
Optimization has been successfully applied to design energy systems across various scales, from utility systems at the plant scale \citep{papouliasStructuralOptimizationApproach1983, vollAutomatedSuperstructurebasedSynthesis2013, bahlRigorousSynthesisEnergy2018, baumgartnerRiSES4RigorousSynthesis2019} to energy systems for districts \citep{bunningBidirectionalLowTemperature2018, schutzOptimalDesignDecentralized2018, teichgraeberClusteringMethodsFind2019} up to power systems on islands \citep{maFeasibilityStudyStandalone2014a, gilsCarbonNeutralArchipelago2017, baroneSupportingSustainableEnergy2021} and on the (inter)-national scale \citep{kannanLongTermElectricityDispatch2013, sialaSustainableEuropeanEnergy2019, reinertEnvironmentalImpactsFuture2020}.
For a review of modeling tools for renewable energy systems, we refer to \cite{ringkjobReviewModellingTools2018}.

Modeling the volatility of VRES leads to uncertainties in the optimization problem parameters.
Two popular approaches to account for uncertain parameters are stochastic programming \citep{dantzigLinearProgrammingUncertainty1955} and robust optimization \citep{campoRobustModelPredictive1987}.
Stochastic programming requires knowledge of the probability distribution of the uncertain parameters and optimizes the expected value of the objective function.
Robust optimization, on the other hand, guarantees the feasibility of an optimal solution for all parameters within a predefined uncertainty set and does not require the probability distribution.
For detailed information on stochastic programming and robust optimization, we refer to \cite{birgeIntroductionStochasticProgramming2011} and \cite{ben-talRobustOptimization2009}, respectively.

To capture the time-varying behavior of VRES, the probability distribution can be approximated discretely by including historical time-series data in the optimization as scenarios in a two-stage stochastic programming \citep{dantzigLinearProgrammingUncertainty1955} formulation of the design problem \citep{teichgraeberTimeseriesAggregationOptimization2022}.
A sufficiently high temporal resolution is needed to properly account for the time-varying behavior of renewable electricity production \citep{ponceletImpactLevelTemporal2016, kelesMeetingModelingNeeds2017, ringkjobReviewModellingTools2018}.
Excessive use of historical data and high temporal resolutions, however, increase the number of constraints and variables of the optimization problem, thus leading to computationally intractable formulations \citep{pfenningerEnergySystemsModeling2014,hoffmannReviewTimeSeries2020}.
Therefore, tradeoffs between computational tractability and the degree of detail in which the system is modeled are necessary.

Different dimensions of model complexity reduction can be leveraged to facilitate computational tractability.
First, model complexity can be reduced by spatial aggregation.
A prominent and extreme example is the copper plate assumption, where the spatial distribution is completely ignored, and transmission losses are neglected \citep{hessRepresentingNodeinternalTransmission2018}.
Furthermore, temporal model complexity can be reduced by using only a representative subset of the historical time-series data, e.g., a collection of representative days \citep{chapaloglouDatainformedScenarioGeneration2022}.
A popular method to determine representative scenarios is clustering \citep{kotzurImpactDifferentTime2018, teichgraeberClusteringMethodsFind2019}.
However, aggregation methods may neglect extreme scenarios, which may impact the operational feasibility of the optimal system design.

Heuristic methods have been introduced to incorporate extreme scenarios during the design.
The heuristics identify and add extreme scenarios to the optimization problem before performing the design optimization.
For example, \cite{dominguez-munozSelectionTypicalDemand2011} exclude days with peak demands from the aggregation and add these to the scenario data of the problem directly.
However, by incorporating individual peak days for the different energy demand forms, e.g., heat, cooling, and electricity, interactions between these demand forms may get lost.
Furthermore, when multiple uncertain quantities are considered, e.g., solar electricity production and electricity demand, the worst-case realization might not correspond to scenarios where the quantities are at their extreme values but rather a scenario where the interaction of their values creates the largest supply and demand mismatch.

Furthermore, the worst-case scenario depends on the energy system design; for example, the installed capacities of solar PV and wind turbines in an energy system influence whether scenarios with low wind speeds or low solar irradiance are critical.
\emph{A priori} extreme period selection cannot consider this influence of the design on the extreme periods.

To mitigate this problem, extreme periods specific to a given (preliminary) design can be identified by optimizing the system \emph{operation}, i.e., by determining whether the considered design can satisfy energy demands for all historical scenarios.
Scenarios in which operations are found to be infeasible are then added as extreme periods to the design problem formulation, and subsequently, another design optimization is performed. 
In the context of energy system design, \cite{bahlTimeseriesAggregationSynthesis2016} introduced such an iterative heuristic in which they repeatedly solve a design problem with approximated operational costs and add virtual time steps for which the operational problem becomes infeasible until eventually a design is obtained that can supply all historical scenarios.
The solution approach that \cite{bahlTimeseriesAggregationSynthesis2016} use to identify worst-case scenarios from a set of finite  cardinality is very similar to the algorithm proposed by \cite{blankenshipInfinitelyConstrainedOptimization1976} to solve semi-infinite programs (SIPs).
\cite{teichgraeberExtremeEventsTime2020} developed this concept further and referred to it as the optimization-based feasibility time-step heuristic.
We will refer to this approach as the \emph{feasibility time-step heuristic} for brevity.
However, as we will show, considering only the historical data may not be enough to identify robust designs.
Specifically, if the operational problem is nonconvex, the identified system design may not be feasible for scenarios that lie ``between'' historical scenarios, i.e., that constitute convex combinations of historical data.

Addressing the aforementioned challenges, we propose a rigorous and robust energy system design (RESD) approach that identifies worst-case scenarios from a predefined space of possible uncertainty realizations \emph{during} the optimization.
The RESD approach bears similarity to the approach for optimal process design under uncertainty by \cite{halemaneOptimalProcessDesign1983}.
We formulate the energy system design problem as a two-stage stochastic program, where we approximate the objective by discretization but treat the constraints continuously by formulating them as semi-infinite constraints \citep{charnesDualityHaarPrograms1962, grossmannOptimumDesignChemical1978, hettichSemiInfiniteProgrammingTheory1993, djelassiRecentAdvancesNonconvex2021}.
The discretized approximation of the objective is obtained by approximating the operational costs using a small number of representative scenarios that we determine using clustering.
To address the neglect of extreme scenarios by clustering, the semi-infinite constraints enforce feasible operation for all scenarios within a predefined uncertainty space and guarantee the robustness of the identified design towards the volatility introduced by the VRES and the energy demand.
Semi-infinite programs often arise in robust design problems \citep{ben-talRobustOptimizationMethodology2002}; for reviews of semi-infinite programming, we refer to \cite{hettichSemiInfiniteProgrammingTheory1993, ruckmannSemiInfiniteProgramming1998, guerravazquezGeneralizedSemiinfiniteProgramming2008, steinHowSolveSemiinfinite2012}.

The key difference between established methods in energy system design, such as the feasibility time-step heuristic, and the RESD approach is the ability of the latter to consider continuous uncertainty spaces, i.e., allow the sets of possible uncertainty realizations to have infinite cardinality.
In contrast, the feasibility time-step heuristic is constrained to finite cardinality sets.
The space of uncertainty realizations can be defined arbitrarily; a straightforward example of such a set is a box-constrained set.
In the present work, we define the uncertainty space by means of the convex hull around the historical data, i.e., we consider convex combinations of historical data as possible uncertainty realizations.

The optimization problem corresponding to the RESD approach is a three-level hierarchical program that is computationally challenging to solve.
If all historical data points are contained within the bounded part of the search space, the optimization problem corresponding to the RESD approach is a restriction of the problem considered in the feasibility time-step heuristic \citep{bahlTimeseriesAggregationSynthesis2016, teichgraeberExtremeEventsTime2020}.
A solution to the restricted problem is robust for more scenarios than the original problem but may have a worse objective function value.

In previous work, we demonstrated that principal component analysis (PCA) is well suited to reduce the dimensionality of energy time-series data \citep{cramerPrincipalComponentDensity2022}.
To facilitate computational tractability, we use PCA \citep{pearsonLinesPlanesClosest1901} to reduce the dimensionality of the data in the RESD approach.
Furthermore, we utilize a lifting approach to improve computational tractability for design problems where the problem of determining the optimal operational strategy is convex. Specifically, we leverage that the lower-level problem is a convex optimization problem to reformulate the problem as a SIP using the approach proposed by \cite{diehlLiftingMethodGeneralized2013}.

The remainder of this work is structured as follows: Section \ref{Methods} introduces the problem structure associated with the RESD approach and explains how the historical data is incorporated into the optimization problem.
In Section \ref{MILP_Case_Study}, we use an illustrative mixed-integer linear problem (MILP) example to demonstrate the shortcomings of the feasibility time-step heuristic in identifying relevant extreme scenarios and demonstrate the robustness of the RESD approach.
Furthermore, special cases in which the feasibility time-step heuristic \citep{bahlTimeseriesAggregationSynthesis2016, teichgraeberExtremeEventsTime2020} yields identical results to the RESD approach are pointed out.
Finally, the lifting approach based on \cite{diehlLiftingMethodGeneralized2013} is introduced to improve computational performance for problems with convex operational behavior.
By means of a case study on the La Palma energy system, Section \ref{Case_Study} then analyzes the accuracy of the designs obtained by the RESD approach and the influences of the dimensionality reduction and the lifting approach on the computational performance.
Finally, Section \ref{Conclusion} summarizes our work.

\section{Robust Energy System Design Approach}\label{Methods}\label{RESD}
Robust energy system designs are feasible for all considered uncertainty realizations, i.e., they are able to supply the energy demand while satisfying all operational constraints.
However, they are not robust with respect to the objective, i.e., they do not give a guaranteed bound on the operational costs.
The RESD approach relies on two concepts to identify robust designs that are cost-optimal: (i) the operational costs are estimated using a small number of representative scenarios with associated probabilities of occurrence \citep{chapaloglouDatainformedScenarioGeneration2022}, (ii) an embedded optimization problem ensures the feasibility for all possible uncertainty realizations within a predefined uncertainty set.

\begin{align*}\tag{PS}\label{problem_structure}
    \underset{\mathbf{x}, \mathbf{z}_s}{\min} & \quad Investment \ costs(\mathbf{x}) + \sum_{s \in S} Operational \ costs(\mathbf{z}_s) ,&\\
    \text{s.t.} & \quad \mathbf{g}_{en}(\mathbf{x}, \mathbf{y_s}, \mathbf{z}_s) \le \mathbf{0} \ \forall s \in \mathcal{S} \ (Energy \ system \ model),&\\ 
    & \quad \underset{\mathbf{y}}{\max}\underset{\mathbf{z}}{\min} \quad Energy \ gap (\mathbf{x}, \mathbf{y}, \mathbf{z}) \le 0 ,&\\
    & \quad \text{s.t.} \quad \mathbf{g}_{en}(\mathbf{x}, \mathbf{y}, \mathbf{z}) \le \mathbf{0} \ (Energy \ system \ model),\\
    & \quad \hskip2.3em\relax \mathbf{g}_{y}(\mathbf{y}) \le \mathbf{0} \ (Uncertainty \ bounds),
\end{align*}

~\eqref{problem_structure} shows the general structure of the RESD problem.
The upper-level problem is the deterministic equivalent of a two-stage stochastic program \citep{birgeIntroductionStochasticProgramming2011} with a discrete approximation of the scenarios by fixed uncertainty realizations $\mathbf{y_s}$ and determines the optimal design decisions $\mathbf{x}$, such as the installed capacities of renewable generators, to achieve a cost-optimal design with respect to the expected value of the operational costs, which is calculated by estimating operational costs with operational decisions $\mathbf{z}_s$.
While we assume that the investment costs are not subject to uncertainty, the operational costs depend on future uncertainty realizations, e.g., on the solar irradiance; thus, we estimate their expected value by means of representative scenarios $\mathcal{S}$ \citep{chapaloglouDatainformedScenarioGeneration2022}.

Furthermore, $\mathbf{g}_{en}(\cdot, \cdot, \cdot)$ are the constraints modeling the energy system, which occur both in the upper- and lower-level problem to ensure the feasibility of the operational decisions associated with the representative scenarios and the extreme scenarios.
\begin{align*}\tag{MLP}\label{medial_level}
    \underset{\mathbf{y}}{\max}\underset{\mathbf{z}}{\min} & \quad Energy \ gap (\mathbf{x}, \mathbf{y}, \mathbf{z}) \le 0 ,&\\
        \text{s.t.} & \quad \mathbf{g}_{en}(\mathbf{x}, \mathbf{y}, \mathbf{z}) \le \mathbf{0} \ (Energy \ system \ model),\\
        & \quad \mathbf{g}_{y}(\mathbf{y}) \le \mathbf{0} \ (Uncertainty \ bounds),
\end{align*}
To ensure robustness, an embedded optimization problem, which we refer to as the medial-level problem \eqref{medial_level}, guarantees that the energy supply gap is nonpositive for all possible uncertainty realizations $\mathbf{y}$, i.e., the energy supply meets or exceeds the energy demand.
This is achieved by identifying the worst-case uncertainty realization for the superordinate design within the predefined uncertainty set defined by $\mathbf{g}_y(\mathbf{y}) \leq \mathbf{0}$.
Reaction to the uncertainty realizations, i.e., recourse actions \citep{birgeIntroductionStochasticProgramming2011}, is considered through operational decisions $\mathbf{z}$, e.g., battery charge and discharge decisions.

\subsection{Processing of Historical Data}\label{sec:data_processing}
The RESD approach requires representative scenarios for operational cost estimation and uncertainty bounds that delimit the search space for the identification of worst-case uncertainty realizations.
We choose to obtain both from historical time-series data.
\begin{figure}
    \centering
    \includegraphics[width=1\linewidth]{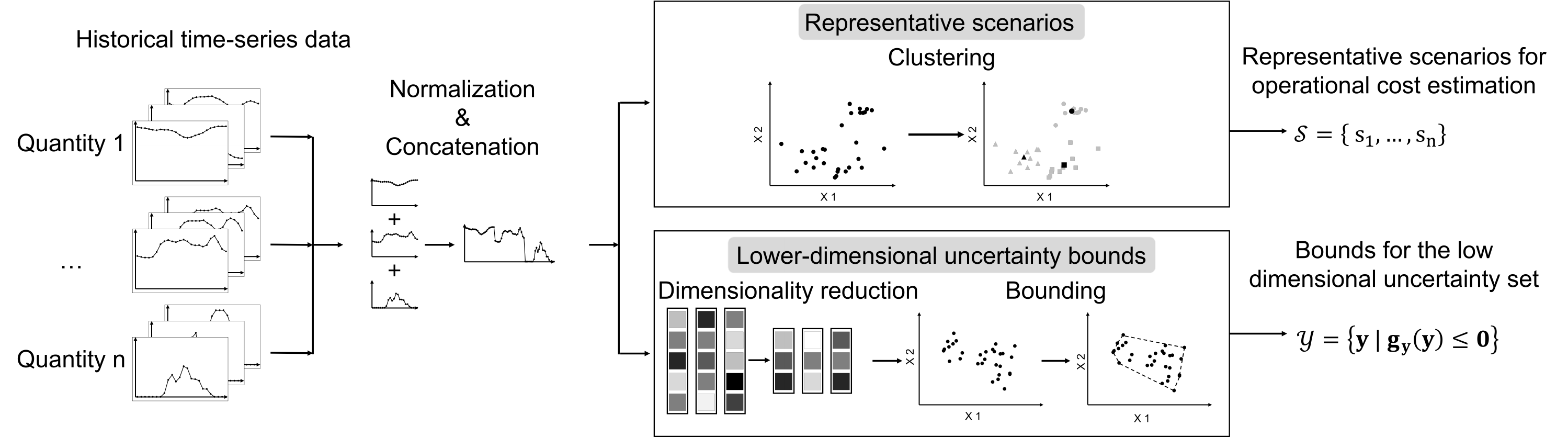}
    \caption{The processing of historical time-series data to obtain the representative scenarios and uncertainty bounds for the RESD problem: Time-series data with the length of the representative period are collected for different quantities, normalized, and concatenated such that a single vector is obtained for each representative period in the historical data. Representative scenarios are obtained by clustering the concatenated and normalized data according to the framework by \cite{teichgraeberClusteringMethodsFind2019}. Uncertainty bounds are obtained by first reducing the dimensionality of the historical data and then determining bounds in the lower-dimensional space.}
    \label{fig:Preprocessing}
\end{figure}
Figure \ref{fig:Preprocessing} shows the corresponding data processing step.
For each quantity of interest, e.g., solar irradiance, the historical time-series data is split into the desired representative period length, e.g., a representative day.

First, representative scenarios are chosen according to the framework by \cite{teichgraeberClusteringMethodsFind2019}, where normalization helps to ensure equal weighting of different quantities during clustering.
The normalized time-series data for each attribute, e.g., solar irradiance and electricity demand, are then concatenated, and representative scenarios are selected by clustering.
Clustering on the concatenated quantities is necessary to maintain the temporal coherence between different quantities.
The choice of the employed normalization technique, clustering method, and parameters is generally problem-specific; we refer to \cite{teichgraeberClusteringMethodsFind2019} for guidance.

Second, the bounds for the uncertainty set need to be determined.
In previous work, we showed that time-series scenarios for renewable energy sources lie on linear lower-dimensional manifolds \citep{cramerPrincipalComponentDensity2022}.
Accordingly, we introduce a dimensionality reduction step to facilitate better computational tractability of the RESD problem.
The dimensionality reduction shifts the worst-case search to a lower-dimensional latent space, reducing the number of variables in the embedded optimization problem.
We choose principal component analysis (PCA) \citep{pearsonLinesPlanesClosest1901} as a linear dimensionality reduction technique.
In principle, nonlinear dimensionality reduction methods could be used.
However, the corresponding reverse transformation is embedded in the optimization problem.
Hence, using a nonlinear reverse transformation gives rise to a nonlinear and, thus, typically much more challenging optimization problem.
Depending on the truncation error, the identified extreme scenarios can be more or less extreme in the original space, meaning that the identified designs could be over- or under-conservative.

After reducing the dimensionality of the historical data by a transformation into the principal component space, we need to determine bounds on the lower-dimensional uncertainty realizations.
The use of nonlinear bounding techniques would give rise to nonlinear optimization problems.
Instead, we choose the convex hull to bound the uncertainty realizations using linear equations.
The convex hull can be defined as a convex combination of its vertices or as points that lie inside the combination of half-spaces, i.e., the faces of the convex hull \citep{avisHowGoodAre1997}; we refer to Section $2.2$ of the supplementary material for more detailed information.
Determining the convex hull has exponential time complexity with respect to the dimensionality of the data \citep{chazelleOptimalConvexHull1993a}, so that for a larger number of principal components the time to compute the convex hull can become a bottleneck, even in comparison to the optimization runtime.

\subsection{The RESD Approach Problem Formulation}\label{RESD Problem}
In the following, constraints are hidden in the sets from which the optimization variables can be chosen for better readability.
The \eqref{robust_design_problem} problem  has the following form:
\begin{alignat*}{4}\tag{RESD}\label{robust_design_problem}
    &\underset{\mathbf{x} \in \mathcal{X}, \mathbf{z}_s \in \mathcal{Z}_s(\mathbf{x}, \mathbf{y_s})}{\min} &\qquad & f(\mathbf{x}) + \sum_{s \in \mathcal{S}} f_{o}(\mathbf{z}_{s}),&\\
    &\text{s.t.} &&\underset{\mathbf{y}\in \mathcal{Y}(\mathbf{x})}{\max}\underset{\mathbf{z}\in \mathcal{Z}(\mathbf{x}, \mathbf{y})}{\min} \underset{t \in \mathcal{T}}{\max} \  e_t(\mathbf{x},\mathbf{y},\mathbf{z}) \le 0
\end{alignat*}
\eqref{robust_design_problem} is a variant of the problem proposed by \cite{halemaneOptimalProcessDesign1983} for optimal process design under uncertainty.
It includes a cost function that depends on both design costs $f(\mathbf{x})$ and approximated operational costs $\sum_{s \in \mathcal{S}} f_{o}(\mathbf{z}_{s})$.
Furthermore, a constraint enforces feasibility in the worst-case: $\underset{\mathbf{y}\in \mathcal{Y}(\mathbf{x})}{\max}\underset{\mathbf{z}\in \mathcal{Z}(\mathbf{x}, \mathbf{y})}{\min} \underset{t \in \mathcal{T}}{\max} \ e_t(\mathbf{x},\mathbf{y},\mathbf{z}) \le 0$.
Here, $\mathbf{x}$ are the design decision variables, e.g., installed PV capacity, and $\mathbf{z}_s$ are the operational decision variables, e.g., battery charging or discharging rates.
The subscript $s$ is a scenario index and indicates that $\mathbf{z}_s$ are the operational decisions associated with the respective scenario, while the fixed parameters $\mathbf{y_s}$ are the uncertainty realizations for scenario $s$.

The robustness of the design is enforced by an embedded optimization problem, which requires that the energy supply gap $e_t(\mathbf{x},\mathbf{y},\mathbf{z})$, i.e., the difference between demand and supply of some energy form, e.g., electricity, has to be nonpositive, i.e., the supply must be equal to or exceed the demand, at all time steps $t$ in the set of considered time steps $\mathcal{T}$.
In the case of multiple demands, e.g., electricity and heat ($n=2$), $e_t(\mathbf{x},\mathbf{y},\mathbf{z})$ will be of the form 
\begin{equation*}
    e_t(\mathbf{x},\mathbf{y},\mathbf{z}) = \max\left(e_{t, 1}(\mathbf{x},\mathbf{y},\mathbf{z}), ..., e_{t, n}(\mathbf{x},\mathbf{y},\mathbf{z}) \right),
\end{equation*}
with $n$ referring to the number of energy forms considered.
The inequality constraint requiring the energy supply gap to be nonpositive has to hold for all possible uncertainty realizations $\mathbf{y}$, e.g., capacity factors for PV, while the operational variables $\mathbf{z}$ can be adjusted to help meet the demand.

The set $\mathcal{X} = \{\mathbf{x} \in \mathbb{R}^{n_x} | \ \mathbf{h}_x(\mathbf{x})=\mathbf{0}, \mathbf{g}_x(\mathbf{x}) \le \mathbf{0}\}$ contains all feasible realizations of the design decisions and the set $\mathcal{Y}(\mathbf{x}) = \{\mathbf{y} \in \mathbb{R}^{n_y} | \ \mathbf{h}_y(\mathbf{x}, \mathbf{y})=\mathbf{0}, \mathbf{g}_y(\mathbf{x}, \mathbf{y}) \le \mathbf{0}\}$ describes the feasible set of the uncertainty realizations.
Sets $\mathcal{Z}_s(\mathbf{x}) = \{\mathbf{z}_s \in \mathbb{R}^{n_{z_s}} | \ \mathbf{h}_{z_s}(\mathbf{x}, \mathbf{y_s}, \mathbf{z}_s)=\mathbf{0}, \mathbf{g}_{z_s}(\mathbf{x}, \mathbf{y_s}, \mathbf{z}_s) \le \mathbf{0}\}$ and $\mathcal{Z}(\mathbf{x}, \mathbf{y}) = \{\mathbf{z} \in \mathbb{R}^{n_z} | \ \mathbf{h}_{z}(\mathbf{x}, \mathbf{y}, \mathbf{z})=\mathbf{0}, \mathbf{g}_{z}(\mathbf{x}, \mathbf{y}, \mathbf{z}) \le \mathbf{0}\}$ describe feasible operational states and are bounded by equations describing the technical behavior of the system, e.g., equations describing the charging and discharging capabilities of the battery, which additionally depend on the capacities of the installed components $\mathbf{x}$.
The dimensions of the vectors are denoted by $n_x, n_{z_s}, n_y, n_z$, respectively.

The embedded optimization problem in \eqref{robust_design_problem} can be replaced by a semi-infinite existence constraint \citep{djelassiDiskretisierungsbasierteAlgorithmenFur2020}: 
\begin{align*}
    \underset{\mathbf{y}\in \mathcal{Y}(\mathbf{x})}{\max}\underset{\mathbf{z}\in \mathcal{Z}(\mathbf{x}, \mathbf{y})}{\min} \underset{t \in \mathcal{T}}{\max} \  e_t(\mathbf{x},\mathbf{y},\mathbf{z}) \le 0 \Longleftrightarrow \forall \mathbf{y}\in \mathcal{Y}(\mathbf{x}) [ \exists \ \mathbf{z}\in \mathcal{Z}(\mathbf{x}, \mathbf{y}) : \underset{t \in \mathcal{T}}{\max} \  e_t(\mathbf{x},\mathbf{y},\mathbf{z}) \le 0 ]
\end{align*}
If coupling inequality constraints, i.e., constraints of the form $\mathbf{g}_y(\mathbf{x}, \mathbf{y}) \leq \mathbf{0}$ that involve the design decisions $\mathbf{x}$ and the uncertainty realizations $\mathbf{y}$, are present, the resulting problem is an existence-constrained generalized semi-infinite program (EGSIP) which can be solved by solving its existence-constrained semi-infinite program (ESIP) relaxation under some mild assumptions \citep{mitsosGlobalOptimizationGeneralized2015}.

Note that special care must be taken if coupling equality constraints are present. In the medial-level problem, coupling equality constraints can contain both upper- and medial-level variables: $\mathbf{h}_y(\mathbf{x}, \mathbf{y}) = \mathbf{0}$. In the lower-level problem, coupling equality constraints can contain upper- and/or medial-level and lower-level variables: $\mathbf{h}_z(\mathbf{x}, \mathbf{y}, \mathbf{z}) = \mathbf{0}$.
Problems with coupling equality constraints may require specialized algorithms since the convergence of the algorithms employed in this work is no longer guaranteed.
However, the coupling equality constraints can be eliminated if they can be rearranged as explicit functions that define some dependent optimization variables as a function of the remaining (independent) variables.
For example, consider a medial-level problem with coupling constraints $\mathbf{h}_y(\mathbf{x}, \mathbf{y}) = \mathbf{0} $.
If the uncertainty realizations $\mathbf{y}$ can be split into a dependent part $\tilde{\mathbf{y}}$ and an independent part $\bar{\mathbf{y}}$ and the coupling equality constraints $\mathbf{h}_y(\mathbf{x}, \mathbf{y}) = \mathbf{h}_y(\mathbf{x}, \tilde{\mathbf{y}}, \bar{\mathbf{y}})  = \mathbf{0}$ can be rearranged into an explicit equation $\tilde{\mathbf{y}} = \tilde{\mathbf{y}}(\mathbf{x}, \bar{\mathbf{y}})$, $\tilde{\mathbf{y}}$
can be replaced by $\tilde{\mathbf{y}}(\mathbf{x}, \bar{\mathbf{y}})$ everywhere, thus eliminating the dependent variables and removing the coupling equality constraints.
\cite{halemaneOptimalProcessDesign1983} use this technique to replace the state variables in the operational stage by a function of the design variables, uncertainties, and operational variables.
In the case that the coupling constraints cannot be rearranged as an explicit function, methods for solving SIPs containing implicit functions, such as the specialized algorithms proposed by \cite{djelassiDiscretizationbasedAlgorithmsGeneralized2019} and \cite{stuberSemiInfiniteOptimizationImplicit2015a} may be used.

To solve the \eqref{robust_design_problem} problem, we first reformulate it as an EGSIP and then form the ESIP relaxation, which we solve using the algorithm introduced by \cite{djelassiGlobalSolutionSemiinfinite2021} and implemented in the libDIPS software package \citep{zinglerLibDIPSDiscretizationBasedSemiInfinite2023}.
To this end, we transform the embedded optimization problem, i.e., the $\underset{\mathbf{y}\in \mathcal{Y}(\mathbf{x})}{\max}\underset{\mathbf{z}\in \mathcal{Z}(\mathbf{x}, \mathbf{y})}{\min} \underset{t \in \mathcal{T}}{\max} \  e_t(\mathbf{x},\mathbf{y},\mathbf{z})$ problem, into a bi-level problem by introducing an auxiliary variable $e_{epi}$ to the minimization problem, which moves the innermost maximization problem into the constraints of the superordinate minimization problem.
Since $\mathcal{T}$ is a finite set, we replace $\underset{t \in \mathcal{T}}{\max} \ e_t(\mathbf{x},\mathbf{y},\mathbf{z}) - e_{epi} \le 0$ by one equation for each time step $t$ and obtain the bi-level problem:
\begin{alignat*}{3}
 &\underset{\mathbf{y}\in \mathcal{Y}(\mathbf{x})}{\max}\underset{\mathbf{z}\in \mathcal{Z}(\mathbf{x}, \mathbf{y}), e_{epi}}{\min} &\qquad  & e_{epi} ,\\
    &\text{s.t.} &&  \ e_t(\mathbf{x},\mathbf{y},\mathbf{z}) - e_{epi} \le 0 &  \ \forall t \in \mathcal{T}
\end{alignat*}
By including the auxiliary variable $e_{epi}$ and the corresponding constraints in the set of lower-level variables $\mathcal{Z}_{epi}(\mathbf{x}, \mathbf{y}) = \{\mathbf{z} \in \mathbb{R}^{n_z}, e_{epi} \in \mathbb{R}| \ \mathbf{h}_{z}(\mathbf{x}, \mathbf{y}, \mathbf{z})=\mathbf{0} \wedge \mathbf{g}_{z}(\mathbf{x}, \mathbf{y}, \mathbf{z}) \le \mathbf{0} \wedge e_t(\mathbf{x},\mathbf{y}, \mathbf{z}) - e_{epi} \le 0 \ \forall t \in \mathcal{T} \}$, we can write the embedded optimization problem as $\underset{\mathbf{y}\in \mathcal{Y}(\mathbf{x})}{\max}\underset{\mathbf{z}\in \mathcal{Z}_{epi}(\mathbf{x}, \mathbf{y})}{\min} e_{epi}$, and thus obtain for the full design problem:
\begin{alignat*}{4}
    &\underset{\mathbf{x} \in \mathcal{X}, \mathbf{z}_s \in \mathcal{Z}_s(\mathbf{x}, \mathbf{y_s})}{\min} &\qquad &f(\mathbf{x}) + \sum_{s \in \mathcal{S}} f_{o}(\mathbf{z}_{s}),&\\
    &\text{s.t.} &&\underset{\mathbf{y}\in \mathcal{Y}(\mathbf{x})}{\max}\underset{\mathbf{z}\in \mathcal{Z}_{epi}(\mathbf{x}, \mathbf{y})}{\min} & e_{epi} \le 0
\end{alignat*}
This problem can now be reformulated as an EGSIP \citep{djelassiGlobalSolutionSemiinfinite2021}:
\begin{alignat*}{4}\tag{EGSIP}\label{EGSIP}
    &\underset{\mathbf{x} \in \mathcal{X}, \mathbf{z}_s \in \mathcal{Z}_s(\mathbf{x}, \mathbf{y_s})}{\min} &\qquad &f(\mathbf{x}) + \sum_{s \in \mathcal{S}} f_{o}(\mathbf{z}_{s}),&\\
    &\text{s.t.} &&\forall \mathbf{y}\in \mathcal{Y}(\mathbf{x}) [ \exists \ \mathbf{z}\in \mathcal{Z}_{epi}(\mathbf{x}, \mathbf{y}) : e_{epi} \le 0 ]&
\end{alignat*}

Generally, if the lower-level feasible set $\mathcal{Z}(\mathbf{x}, \mathbf{y})$ of the \eqref{robust_design_problem} problem is defined by constraints of the form $\mathbf{g}_z(\mathbf{x}, \mathbf{y}, \mathbf{z}) \le \mathbf{0}$, a reformulation of the problem is required to ensure applicability of the ESIP algorithm \citep{djelassiDiskretisierungsbasierteAlgorithmenFur2020}.
This reformulation moves the constraints $\mathbf{g}_z(\mathbf{x}, \mathbf{y}, \mathbf{z}) \le \mathbf{0}$ into the objective function, which enables the medial level to identify $\mathbf{y}$ that make the lower-level problem infeasible.
However, in problem \eqref{EGSIP}, the epigraph reformulation constraints $e_t(\mathbf{x},\mathbf{y},\mathbf{z}) - e_{epi} \le 0  \ \forall t \in \mathcal{T}$ are always feasible if the bounds on $e_{epi}$ are sufficiently large, which we can always guarantee by choosing the bounds appropriately.
The guaranteed feasibility of the epigraph constraints allows us to omit the reformulation step for the epigraph reformulation constraints.
Similarly, it is not necessary to move constraints of the form $\mathbf{g}_z(\mathbf{x}, \mathbf{z}) \le \mathbf{0}$ into the objective function since $\mathbf{y}$ cannot cause infeasibility of  $\mathbf{g}_z(\mathbf{x}, \mathbf{z})$.
Note that uncertainties typically occurring in energy system design problems are uncertain demand and uncertain generation by renewable energy sources, which enter the problem only in the energy balance, i.e., the objective function of the medial-level problem.

The next step is to form the ESIP relaxation.
To this end, we first write problem \eqref{EGSIP} equivalently \citep{djelassiDiskretisierungsbasierteAlgorithmenFur2020} as 
\begin{alignat*}{4}
    &\underset{\mathbf{x} \in \mathcal{X}, \mathbf{z}_s \in \mathcal{Z}_s(\mathbf{x}, \mathbf{y_s})}{\min} &\qquad &f(\mathbf{x}) + \sum_{s \in \mathcal{S}} f_{o}(\mathbf{z}_{s}),&\\
    &\text{s.t.} &&\forall \mathbf{y}\in \mathcal{Y}_{ref} [ \exists \ \mathbf{z}\in \mathcal{Z}_{epi}(\mathbf{x}, \mathbf{y}) : e_{epi} \le 0 \lor \mathbf{g}_y(\mathbf{x}, \mathbf{y}) > \mathbf{0}],&
\end{alignat*}
where $\mathcal{Y}_{ref} = \{\mathbf{y} \in \mathbb{R}^{n_y} | \ \mathbf{h}_y(\mathbf{y})=\mathbf{0}, \mathbf{g}_y(\mathbf{y}) \le \mathbf{0}\}$ is the set of possible uncertainty realizations without the inequality constraints that depend on the upper-level variables $\mathbf{x}$, which instead have been moved into the logical constraint $e_{epi} \le 0 \lor \mathbf{g}_y(\mathbf{x}, \mathbf{y}) > \mathbf{0}$.
This logical constraint describes that either the medial-level problem is feasible and the energy supply gap is nonpositive or the medial-level problem is infeasible, i.e., the set of worst-case uncertainty realizations is empty.
In either case, the upper-level problem is feasible.
Furthermore, we assume here that the coupling equality constraints $\mathbf{h}_y(\mathbf{x}, \mathbf{y})=\mathbf{0}$, if they were present, have been removed using the procedure described at the beginning of this Section.

To make the ESIP algorithms available in libDIPS \citep{zinglerLibDIPSDiscretizationBasedSemiInfinite2023} applicable, the strictness of the inequality $-\mathbf{g}_y(\mathbf{x}, \mathbf{y}) < \mathbf{0}$ is relaxed, resulting in an ESIP problem:
\begin{alignat*}{4}
    &\underset{\mathbf{x} \in \mathcal{X}, \mathbf{z}_s \in \mathcal{Z}_s(\mathbf{x}, \mathbf{y_s})}{\min} &\qquad &f(\mathbf{x}) + \sum_{s \in \mathcal{S}} f_{o}(\mathbf{z}_{s}),&\\
    &\text{s.t.} &&\forall \mathbf{y}\in \mathcal{Y}_{ref} [ \exists \ \mathbf{z}\in \mathcal{Z}_{epi}(\mathbf{x}, \mathbf{y}) : e_{epi} \leq 0 \lor  - \mathbf{g}_y(\mathbf{x}, \mathbf{y}) \le \mathbf{0}],&
\end{alignat*}
which can also be written as:
\begin{alignat*}{4}\tag{ESIP REL}\label{ESIP Rel}
    &\underset{\mathbf{x} \in \mathcal{X}, \mathbf{z}_s \in \mathcal{Z}_s(\mathbf{x}, \mathbf{y_s})}{\min} &\qquad &f(\mathbf{x}) + \sum_{s \in \mathcal{S}} f_{o}(\mathbf{z}_{s}),&\\
    &\text{s.t.} &&\forall \mathbf{y}\in \mathcal{Y}_{ref} [ \exists \ \mathbf{z}\in \mathcal{Z}_{epi}(\mathbf{x}, \mathbf{y}) : \min(e_{epi}, -{g}_{y, j}(\mathbf{x}, \mathbf{y}) \, | \, j \in \{1, ..., n_{g_y}\}) \le 0]&
\end{alignat*}
Note that this relaxation is generally inexact.
However, it is established in SIP literature that for all but degenerate cases, the relaxed problem should lead to the same objective value as the original problem \cite{mitsosGlobalOptimizationGeneralized2015}.

For notational convenience, we define $g_{e}(\mathbf{x}, \mathbf{y}, \mathbf{z}) := \min(e_{epi}, -{g}_{y, j}(\mathbf{x}, \mathbf{y}) \, | \, j \in \{1, ..., n_{g_y}\})$ and write \ref{ESIP Rel} as:
\begin{alignat*}{4}\tag{ESIP}\label{ESIP}
    &\underset{\mathbf{x} \in \mathcal{X}, \mathbf{z}_s \in \mathcal{Z}_s(\mathbf{x}, \mathbf{y_s})}{\min} &\qquad &f(\mathbf{x}) + \sum_{s \in \mathcal{S}} f_{o}(\mathbf{z}_{s}),&\\
    &\text{s.t.} &&\forall \mathbf{y}\in \mathcal{Y}_{ref} [ \exists \ \mathbf{z}\in \mathcal{Z}_{epi}(\mathbf{x}, \mathbf{y}): g_{e}(\mathbf{x}, \mathbf{y}, \mathbf{z}) \le 0]&
\end{alignat*}

The \eqref{ESIP} problem determines design decision variables $\mathbf{x}$ and operational decision variables $\mathbf{z}_s$ such that for all feasible uncertainty realizations $\mathbf{y}$ in the set of possible uncertainty realizations $\mathcal{Y}_{ref}$ there exist operational recourse actions $\mathbf{z}$ so that the constraint $g_{e}(\mathbf{x}, \mathbf{y}, \mathbf{z}) \le 0$ is satisfied, meaning that either $e_{epi} \le 0$ is satisfied or one of the EGSIP constraints is violated, i.e., $-{g}_{y, j}(\mathbf{x}, \mathbf{y}) \le 0$, and the set of uncertainty realizations of the EGSIP $\mathcal{Y}(\mathbf{x})$ is empty.

\section{Limitations of the Feasibility Time-Step Heuristic}\label{MILP_Case_Study}
This section highlights how the feasibility time-step heuristic cannot rigorously provide robust designs for nonconvex problems.
Specifically, we showcase in Section \ref{MILP Example} how the RESD approach can identify a robust design for a nonconvex problem while the feasibility time-step heuristic \citep{bahlTimeseriesAggregationSynthesis2016, teichgraeberExtremeEventsTime2020} fails.
Section \ref{special} shows the cases for which the feasibility time-step heuristic yields robust designs.
Finally, in Section \ref{lifting}, we present an extension to the RESD approach that improves performance for problems with convex lower-level problem.

\subsection{MILP Example}\label{MILP Example}
Ensuring operability for all historical data may not be enough to obtain a robust design if the operational problem is nonconvex.
SIPs with nonconvex lower-level problems are challenging to solve and cannot, in general, be reformulated as finite optimization problems \citep{djelassiRecentAdvancesNonconvex2021}.
The reformulation into a tractable finite optimization problem requires the equations and uncertainty set to fulfill special criteria \citep{bertsimasTheoryApplicationsRobust2011}.
MILPs are nonconvex problems and are extensively used in energy system design \citep{pfenningerEnergySystemsModeling2014}.
Therefore, we use a small MILP example that mimics an energy system to demonstrate the shortcomings of the feasibility time-step heuristic and solve it using the RESD approach.

The MILP problem reads
\begin{alignat*}{4}\tag{MILP RESD}\label{milp_example}
    &\underset{(x_1, x_2) \in [0, 100]^2}{\min} &\qquad & 2x_1 + x_2,&\\
    &\text{s.t.} &&\forall y_1 \in [0, 100][\exists (z_1, z_2, b) \in \mathcal{Z}(x_1, x_2, y_1) : y_1 - z_1 - z_2 \le 0],
\end{alignat*}
with the feasible set
\begin{align*}
    \mathcal{Z}(x_1, x_2, y_1) = \{(z_1, z_2) \in [0, 100]^2, b \in \{0, 1\}\ | \ z_1 - x_1 \le 0 \wedge z_2 - b x_2 \le 0 \wedge\\
    0.2 b x_2 - z_2 \le 0 \wedge y_1 - z_1 - z_2 = 0\}.
\end{align*}

In \eqref{milp_example}, $x_1$ and $x_2$ represent abstract component capacities, $y_1$ an uncertain demand, and $z_1$ and $z_2$ component outputs, while $b$ models on/off decisions for component $2$ to account for the minimal part load restriction of that component.
Compared with Component $1$, Component $2$ is half as expensive to install but has a minimum part load, below which the component has to be shut off.

For simplicity, we assume cost and demand data without reference to real-world data. Specifically, Component 1 is assumed to cost $2$ units, and Component 2 is assumed to cost $1$ unit.
Demand is assumed to be in the range between $0$ and $100$ units and has to be met exactly without curtailment of excess power.

Note that the coupling equality constraint $y_1 - z_1 - z_2 = 0$, which represents the no-curtailment constraint, can be handled by applying the approach described in Section \ref{RESD Problem}, i.e., 
we solve the equality constraint explicitly for $z_2$ and substitute $z_2 = y_1 - z_1$ to obtain the following problem:
\begin{alignat*}{4}\tag{MILP FEAS}\label{milp_example_feas}
    &\underset{(x_1, x_2) \in [0, 100]^2}{\min} &\qquad & 2x_1 + x_2,&\\
    &\text{s.t.} &&\forall y_1 \in [0, 100][\exists (z_1, b) \in \mathcal{Z}_{feas}(x_1, x_2, y_1):0 \le 0],
\end{alignat*}
with the feasible set
\begin{align*}
    \mathcal{Z}_{feas}(x_1, x_2, y_1) = \{z_1 \in [0, 100], b \in \{0, 1\}\ | \ z_1 - x_1 \le 0 \wedge y_1 - z_1 - b x_2 \le 0 \wedge 0.2 b x_2 - y_1 + z_1 \le 0\}.
\end{align*}
When substituting $z_2 = y_1 - z_1$ into the SIP constraint, we obtain the trivial constraint $0 \le 0$, which is always satisfied if the feasible set $\mathcal{Z}_{feas}(x_1, x_2, y_1)$ is not empty.
Hence, \eqref{milp_example_feas} ensures that the lower-level feasible set is not empty.

The lower-level feasible set $\mathcal{Z}_{feas}(x_1, x_2, y_1)$ depends on the uncertain variable $y_1$.
As described in Section \ref{RESD Problem}, we remove this dependency by moving the constraints into the objective function of the medial- and lower-level problem, i.e., 
\begin{alignat*}{4}\tag{MILP ESIP REF}\label{esip_relaxation}
    &\underset{(x_1, x_2) \in [0, 100]^2}{\min} &\qquad & 2x_1 + x_2,&\\
    &\text{s.t.} &&\forall y_1 \in [0, 100][\exists (z_1, b) \in \mathcal{Z}_{ref}(x_1, x_2) : \max\{y_1 - z_1 - b x_2, 0.2 b x_2 - y_1 + z_1\} \le 0],
\end{alignat*}
with the feasible set
\begin{align*}
    \mathcal{Z}_{ref}(x_1) = \{z_1 \in [0, 100], b \in \{0, 1\}\ | \ z_1 - x_1 \le 0\}.
\end{align*}
and solve it using the ESIP algorithm implemented in libDIPS \citep{zinglerLibDIPSDiscretizationBasedSemiInfinite2023}.
The involved sub-problems and the solver settings can be found in Section \si{1} and Table \si{2} of the supplementary material, respectively. 

We obtain an optimal system design with $x_1 = 16.60$ and $x_2 = 83.40$, i.e., the capacity of the more expensive Component 1 is roughly \SI{20}{\percent} of the capacity of Component 2.
The fact that the more expensive Component 1 is built allows the system to satisfy the demand below the minimal part load of Component 2.

If only the vertices of the uncertain demand set $\{0, 100\}$ had been considered, the optimal solution would have been $x_1 = 0$ and $x_2 = 100$.
However, this solution is not robust to certain possible intermediate demands.
For example, the demand $y_1 = 15$ lies below the minimal part load of Component 2, i.e., $0.2 \cdot 100 = 20$. Since Component 1 would have been installed, this demand could not be satisfied.

Note that the nonconvexity introduced by the minimal part load of Component $2$ would not lead to infeasibility without the no-curtailment assumption.
If curtailment were allowed, Component $2$ could be operated at minimal part load for demands smaller than its minimal part load.
Consequently, the minimal part load could be neglected in the medial-level problem, effectively removing the nonconvexity.
Thus, nonconvexities do not necessarily lead to worst-case scenarios that lie in-between historical scenarios.
However, the example highlights that care must be taken when applying the feasibility time-step heuristic \citep{bahlTimeseriesAggregationSynthesis2016, teichgraeberExtremeEventsTime2020} to nonconvex problems.

\subsection{When Considering Historical Data is Sufficient}\label{special}
In the previous section, we highlighted how the feasibility time-step heuristic \citep{bahlTimeseriesAggregationSynthesis2016, teichgraeberExtremeEventsTime2020} may fail for nonconvex problems.
However, there are special cases for which the heuristic yields designs that are identical to the ones found by the RESD approach using the convex hull bounding the full dimensional uncertainty space.
In these special cases, the objective function $g_e(\mathbf{x}, \mathbf{y}, \mathbf{z})$ and the constraints of the lower-level problem need to be such that worst-case scenarios lie at an extreme point of $\mathcal{Y}_{ref}$.
An analogy to this is that the maximum of a convex function lies at an extreme point of its domain.
The identified designs are then robust to all uncertainty realizations in the convex hull of the considered historical data.

Worst-case scenarios can then be identified by solving the embedded \eqref{MAXMIN} problem, i.e.,
\begin{equation}\tag{MAXMIN}\label{MAXMIN}
   \underset{\mathbf{y}\in \mathcal{Y}_{ref}}{\max}\underset{\mathbf{z}\in \mathcal{Z}_{epi}(\mathbf{x}, \mathbf{y})}{\min} g_e(\mathbf{x}, \mathbf{y}, \mathbf{z}),
\end{equation}
using vertex enumeration \citep{halemaneOptimalProcessDesign1983}, i.e., the vertices of the feasible set of uncertainty realizations $\mathcal{Y}_{ref}$ are enumerated, and the operational problem is solved for each vertex.
It is thus sufficient to examine the historical data that constitute the vertices of the convex hull and unnecessary to consider the entire historical data set, as it is done in the feasibility time-step heuristic.

There are two special cases for which \cite{halemaneOptimalProcessDesign1983} and \cite{bialasTwolevelOptimization1982} have proven that the worst-case scenarios lie on the vertices of the feasible set and for which the feasibility time-step heuristic thus leads to robust designs: (i) if the semi-infinite constraint $g_e(\mathbf{x},\mathbf{y},\mathbf{z})$ is jointly convex in $\mathbf{y}$ and $\mathbf{z}$ and all other lower-level constraints are convex in $\mathbf{z}$ \citep{halemaneOptimalProcessDesign1983}, and (ii) if both the medial-level and the lower-level problem of the \eqref{robust_design_problem} problem are linear, i.e., the objectives and constraints are linear and consequently the embedded MAXMIN problem is a bilevel linear program (BLLP) \citep{bialasTwolevelOptimization1982}.

\subsection{Lifting Approach for Convex Lower-Level Problems}\label{lifting}
Finally, we introduce a lifting approach to improve the computational performance of the RESD approach that is applicable if the lower-level problem is convex, i.e., the operational problem is convex.
The motivation for this reformulation is that solving \eqref{ESIP} using the ESIP algorithm proposed by \cite{djelassiGlobalSolutionSemiinfinite2021}, which is based on the \cite{blankenshipInfinitelyConstrainedOptimization1976} algorithm, proved to be computationally slow in our preliminary investigations.
Specifically, we found that the majority of the CPU time to solve the ESIP was spent on the solution of the embedded optimization problem \eqref{MAXMIN}, which is equivalent to the semi-infinite existence constraint of \eqref{ESIP}, i.e.,  $\forall \mathbf{y}\in \mathcal{Y}_{ref} [ \exists \ \mathbf{z}\in \mathcal{Z}_{epi}(\mathbf{x}, \mathbf{y}): g_{e}(\mathbf{x}, \mathbf{y}, \mathbf{z}) \le 0]$.

The lifting approach described in  \cite{diehlLiftingMethodGeneralized2013} allows us to reformulate \eqref{MAXMIN} as a single-level nonlinear program (NLP).
To this end, we first reformulate \eqref{MAXMIN} as a SIP by writing its epigraph reformulation with the auxiliary variable $\phi$ \citep{steinBilevelStrategiesSemiinfinite2003}:
\begin{alignat*}{4}\tag{ES}\label{EmbeddedSIP}
    & \underset{\mathbf{y}\in \mathcal{Y}_{ref}, \phi}{\max}  && \phi,&\\
    & \text{s.t.} \qquad && \phi - g_e(\mathbf{x},\mathbf{y},\mathbf{z}) \le 0 \ \forall \mathbf{z} \in \mathcal{Z}_{epi}(\mathbf{x}, \mathbf{y})&
\end{alignat*}
We then employ the lifting method \citep{diehlLiftingMethodGeneralized2013} to reformulate \eqref{EmbeddedSIP} into the NLP
\begin{alignat*}{2}
    & \underset{\mathbf{y}\in \mathcal{Y}_{ref}, \mathbf{z}\in \mathcal{Z}_{epi}(\mathbf{x}, \mathbf{y}), \phi, \boldsymbol{\lambda}, \boldsymbol{\mu}}{\max}  & \phi ,\\
    & \text{s.t.} & \phi -  \mathcal{L}(\mathbf{x},\mathbf{y},\mathbf{z}, \boldsymbol{\lambda}, \boldsymbol{\mu}) = 0 ,\\
    && \nabla_{\mathbf{z}} \phi - \nabla_{\mathbf{z}} \mathcal{L}(\mathbf{x},\mathbf{y},\mathbf{z}, \boldsymbol{\lambda}, \boldsymbol{\mu}) = \mathbf{0},\\
    && \boldsymbol{\mu}\ge \mathbf{0},\\
    && \boldsymbol{\mu}^\intercal\mathbf{g}_z(\mathbf{x},\mathbf{y},\mathbf{z}) = \mathbf{0},
\end{alignat*}
with Lagrange multipliers $\boldsymbol{\lambda}$ and $\boldsymbol{\mu}$ and 
\begin{align*}
    &\mathcal{L}(\mathbf{x},\mathbf{y},\mathbf{z}, \boldsymbol{\lambda}, \boldsymbol{\mu}) = g_e(\mathbf{x},\mathbf{y},\mathbf{z}) + \boldsymbol{\lambda}^\intercal \mathbf{h}_z(\mathbf{z}) + \boldsymbol{\mu}^\intercal \mathbf{g}_z(\mathbf{x},\mathbf{y},\mathbf{z}).
\end{align*}
Note that $\nabla_{\mathbf{z}} \phi = 0$ and by substituting $\phi$  in the objective function using $ \phi =  \mathcal{L}(\mathbf{x},\mathbf{y},\mathbf{z}, \boldsymbol{\lambda}, \boldsymbol{\mu})$ we end up with the single-level problem:
\begin{alignat*}{2}\tag{NLP}\label{EmbeddedNLP}
    & \underset{\mathbf{y}\in \mathcal{Y}_{ref} \mathbf{z}\in \mathcal{Z}_{epi}(\mathbf{x}, \mathbf{y}), \boldsymbol{\lambda}, \boldsymbol{\mu}}{\max}  & \mathcal{L}(\mathbf{x},\mathbf{y},\mathbf{z}, \boldsymbol{\lambda}, \boldsymbol{\mu}) ,\\
    & \text{s.t.} &\nabla_{\mathbf{z}} \mathcal{L}(\mathbf{x},\mathbf{y},\mathbf{z}, \boldsymbol{\lambda}, \boldsymbol{\mu}) = \mathbf{0} ,\\
    && \boldsymbol{\mu}\ge \mathbf{0} ,\\
    && \boldsymbol{\mu}^\intercal\mathbf{g}_z(\mathbf{x}, \mathbf{y}, \mathbf{z}) = \mathbf{0}
\end{alignat*}

The objective is nonlinear because of the multiplication of the Lagrange multipliers $\boldsymbol{\mu}$ and $\boldsymbol{\lambda}$ with the constraints $\mathbf{g}_z(\mathbf{x},\mathbf{y},\mathbf{z})$ and $\mathbf{h}_z(\mathbf{z})$.
Note that, contrary to \cite{diehlLiftingMethodGeneralized2013}, we added the nonlinear complementarity constraints $\boldsymbol{\mu}^\intercal\mathbf{g}_z(\mathbf{x},\mathbf{y},\mathbf{z}) = \mathbf{0}$ since we found this to help convergence in our case study (cf. Section \ref{Case_Study}).
This is an unexpected result since a major motivation of the lifting approach is to avoid the complementarity constraints, which are numerically poorly behaved  \citep{scheelMathematicalProgramsComplementarity2000}.

Finally, by replacing \eqref{MAXMIN} with its single-level reformulation \eqref{EmbeddedNLP}, we obtain the SIP reformulation of the \eqref{robust_design_problem}, i.e., 
\begin{alignat*}{3}\tag{RESD SIP}\label{RESD SIP}
    &\underset{\mathbf{x} \in \mathcal{X}, \mathbf{z}_s \in \mathcal{Z}_s(\mathbf{x})}{\min} &\qquad &f(\mathbf{x}) + \sum_{s \in \mathcal{S}} f_{o}(\mathbf{z}_{s}),\\
    &\text{s.t.}  && \mathcal{L}(\mathbf{x}, \mathbf{w}) \le 0 \ \forall \mathbf{w}\in \mathcal{W}(\mathbf{x}) ,
\end{alignat*} 
with 
\begin{equation*}
    \begin{split}
    \mathcal{W}(\mathbf{x}) = \{\mathbf{y} \in \mathcal{Y}_{ref}, \mathbf{z} \in \mathcal{Z}_{epi}(\mathbf{x}, \mathbf{y}), \boldsymbol{\lambda}, \boldsymbol{\mu} \ | \\
    \nabla_{\mathbf{z}} \mathcal{L}(\mathbf{x},\mathbf{y},\mathbf{z}, \boldsymbol{\lambda}, \boldsymbol{\mu}) = \mathbf{0}, \ \boldsymbol{\mu}\ge \mathbf{0},\\
    \ \boldsymbol{\mu}^\intercal\mathbf{g}_l(\mathbf{x},\mathbf{y},\mathbf{z}) = \mathbf{0}\}.\\
    \end{split}
\end{equation*}

\eqref{RESD SIP} can be solved using the implementation of the Blankenship \& Falk algorithm \citep{blankenshipInfinitelyConstrainedOptimization1976} available in libDIPS \citep{zinglerLibDIPSDiscretizationBasedSemiInfinite2023}.

\section{La Palma Energy System}\label{Case_Study}
We now analyze how close the lower-dimensional RESD problem approximates the costs of the full-dimensional problem as well as the performance of the RESD approach using the example of determining a robust energy system design for the island of La Palma in the Canary Islands.
La Palma is an example of an isolated energy system with significant decarbonization potential.
Isolated energy systems are defined by a lack of connection to a superordinate energy grid and, consequently, a necessity for self-sufficiency in electricity production.
Due to their maturity and reliability, diesel engines have been the electricity generators of choice for many isolated systems \citep{kennedyOptimalHybridPower2017}.
As of 2021, only \SI{20}{\percent} of gross electricity production in the whole Canary Islands was supplied by renewable sources, highlighting the presently strong reliance on fossil fuels \citep{gobiernodecanariasAnuarioEnergeticoCanarias2023}.
Furthermore, the Canary Islands have a high availability of renewable energy resources, especially solar \citep{meschedeClassificationGlobalIsland2016}, and plans to decarbonize their economy and achieve carbon neutrality by 2040 \citep{gobiernodecanariasConsejeriaTransicionEcologica2023} have been announced.

We consider the following system components: PV units, wind turbines, diesel generators, and battery systems.
We neglect to model the electricity distribution system and the related transmission losses to keep the computational complexity manageable; thus, all system components are directly connected to a single node to supply the island's energy demand.

We use data for wind speed and global irradiance obtained from the Photovoltaic Geographical Information System (PVGIS) by the European Union \citep{europeancommissionJRCPhotovoltaicGeographical2022, huldNewSolarRadiation2012}.
Furthermore, electricity demand data was aggregated from the Spanish electricity system operator \cite{redelectricadeespanaPalmaElectricityDemand2024}.
The data covers the timespan from \si{2013} to \si{2019} and is separated into single-day periods with an hourly resolution.
We use 15 representative scenarios in the approximation of the operational costs, as this number leads to a decent approximation of the load curve (see supplementary materials, Figure \si{2}).
Detailed information about the case study, such as the modeling of the individual components, determination of cost parameters, and the overall problem formulation, is provided in Section \si{2} of the supplementary materials. 

The resulting RESD problem has linear medial- and lower-level problems. Hence, the lifting approach (cf. Section \ref{lifting}), the vertex enumeration (cf. Section \ref{special}), and thus also the feasibility time-step heuristic \citep{bahlTimeseriesAggregationSynthesis2016, teichgraeberExtremeEventsTime2020} are applicable.
We utilize the latter as a reference to solve the full-dimensional problem and compare the solution to those of the RESD approach with different degrees of dimensionality reduction.
We do not compare the computational performance of the feasibility time-step heuristic with that of the RESD approach, as the RESD approach is generally applicable, i.e., for nonconvex problems, whereas the feasibility time-step heuristic yields guaranteed robust designs only in the special cases described in Section \ref{special}.

The full problem formulation and the application of the reformulation steps described in Section \ref{Methods} to transform the embedded optimization problem into a single-level NLP are provided in the supplementary materials in Sections $2.2$ \& $2.4$.
The resulting \eqref{ESIP} and \eqref{RESD SIP} are solved with the Blankenship \& Falk-based ESIP algorithm and the Blankenship \& Falk algorithm, respectively, using the implementation in libDIPS \citep{zinglerLibDIPSDiscretizationBasedSemiInfinite2023} and the Gurobi solver version $11.0$ \citep{gurobi}.
We use a desktop computer with a $4$-core/$4$-thread Intel i5-4570 CPU with \SI{3.2}{\giga\hertz}/\SI{3.6}{\giga\hertz} base/turbo frequency and \SI{16}{\giga\byte} of RAM running Microsoft Windows 10 Enterprise version $10.0.17763$.
Optimization settings for libDIPS and Gurobi deviating from the default values are given in Section \si{3} of the supplementary materials; all $4$ cores were allocated to the numerical experiments.

\subsection{Optimal Robust Design}
\begin{figure}[H]
    \centering
    \includegraphics[width=222pt]{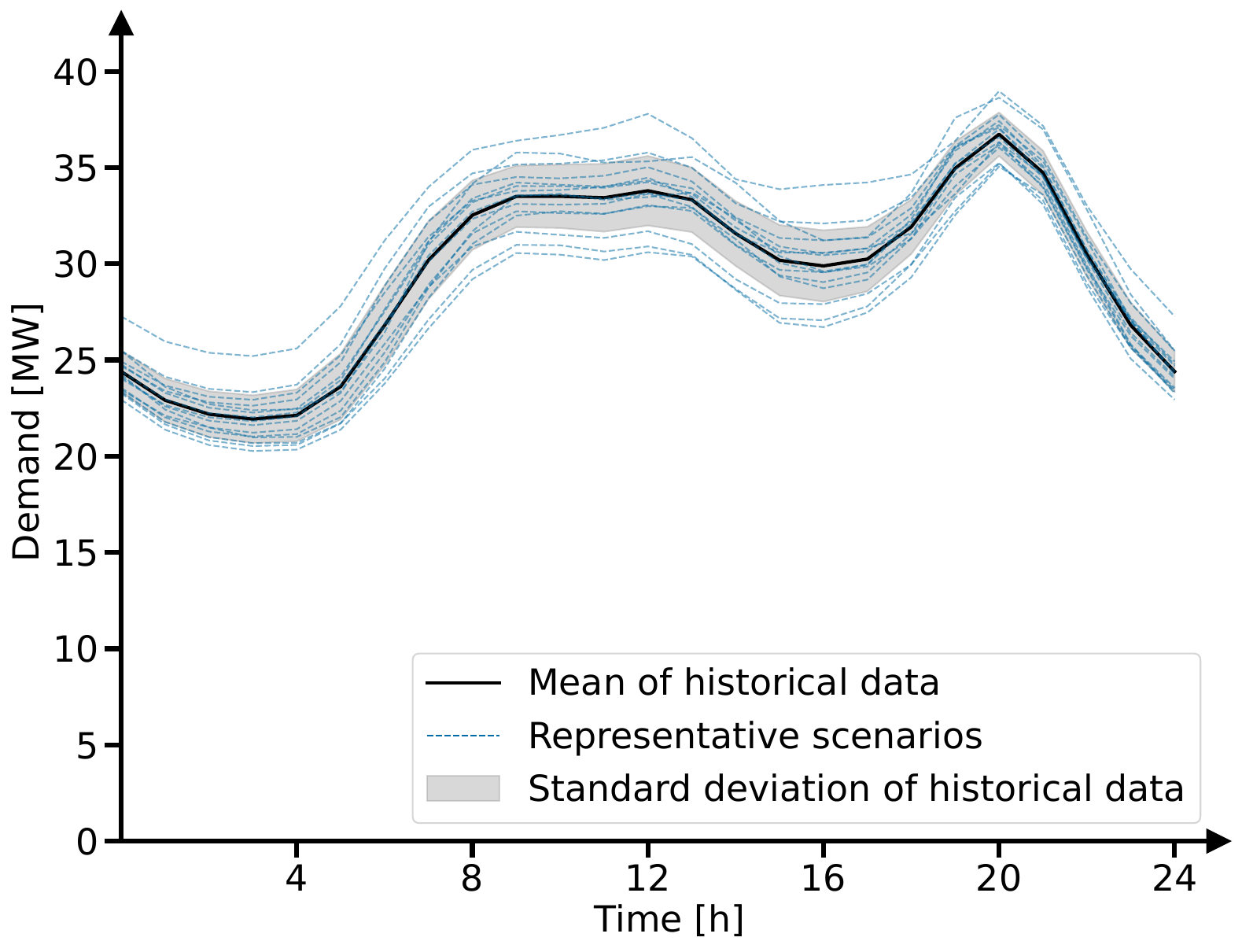}
    \caption{Mean and standard deviation of the daily historical electricity demand on the island of La Palma in 2013-2019. Demand data were obtained from the Spanish electricity distribution system operator \cite{redelectricadeespanaPalmaElectricityDemand2024}. Furthermore, the \si{15} representative scenarios used to compute the operational costs are depicted as dashed lines. There is significant variation in demand over the course of a day, with a wide peak during the morning and a sharp peak in the evening.}
    \label{fig:demand}
\end{figure}

First, we briefly analyze the robust design resulting from the RESD approach and compare it to La Palma's current energy system.
Figure \ref{fig:demand} shows the mean and the standard deviation of the daily electricity demand. Additionally, the \si{15} representative scenarios used to approximate the operational costs are shown.
The demand varies significantly over the course of a day, with a multiple-hour-long peak in the middle of the day between \si{08:00} and \si{12:00} and a sharp peak in the evening around \si{20:00}.
The demand data varies between \SI{0.2}{MW} and \SI{44.8}{MW}, with the low value of \SI{0.2}{MW} being associated with a power outage that occurred on December 11th, 2013.

\begin{figure}[H]
    \centering
    \includegraphics[width=222pt]{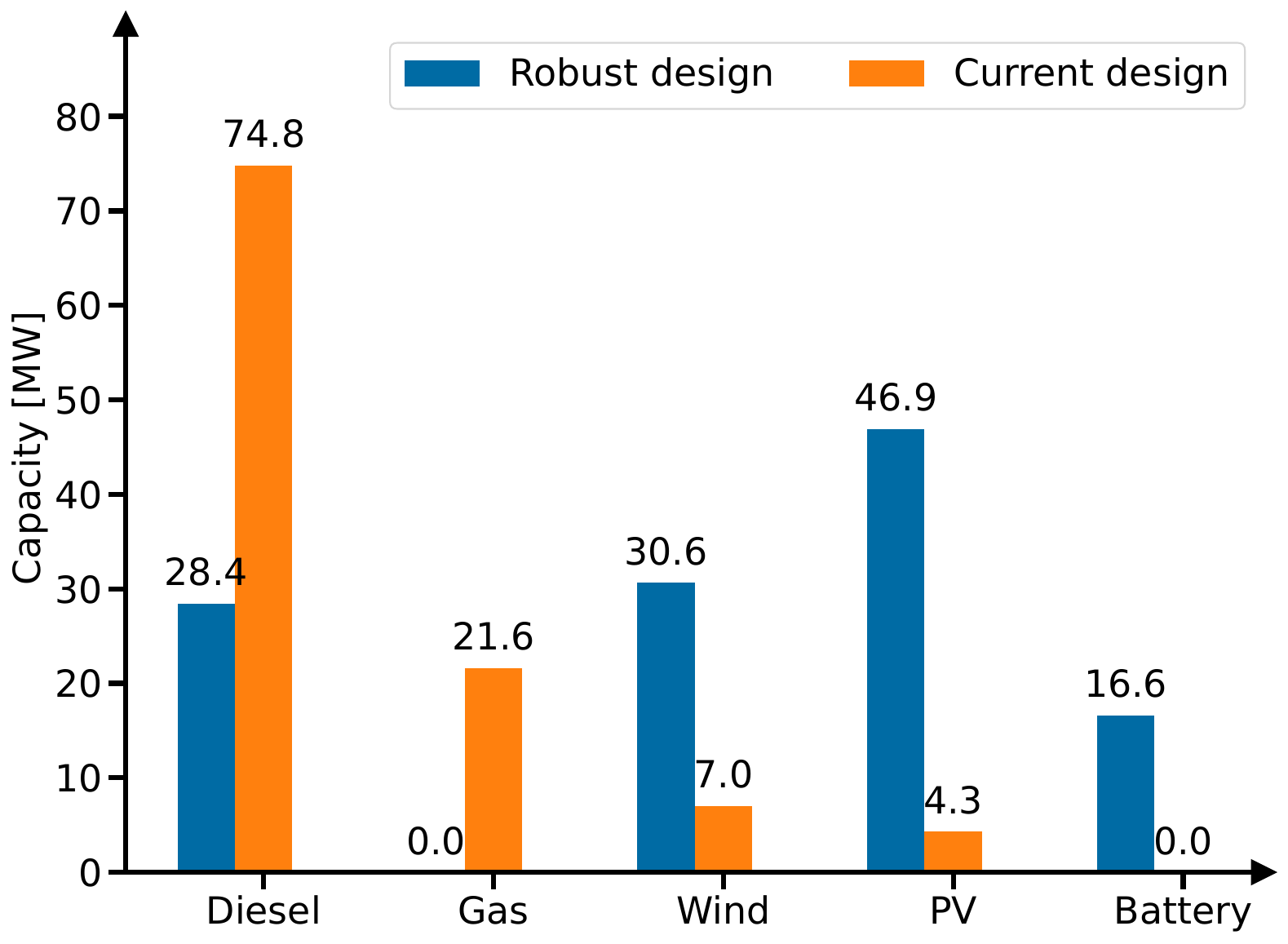}
    \caption{Installed component capacities for the robust design (blue) identified by the RESD approach using \si{16} time steps and \si{9} principal components. For comparison, the currently installed capacities (orange) \citep{gobiernodecanariasAnuarioEnergeticoCanarias2023} are shown. The current design has much higher conventional generation capacity. However, plans to switch to higher renewable electricity generation have already been announced \citep{gobiernodecanariasConsejeriaTransicionEcologica2023}. Note that we did not consider possible component failures in the RESD approach, and hence do not have backup capacities, which are considered in the current design.}
    \label{fig:comparison_capacities}
\end{figure}

Figure \ref{fig:comparison_capacities} shows the currently installed capacities as well as the capacities of the robust design obtained by the RESD approach using \si{16} time steps and \si{9} principal components.
The total annualized costs (TAC) of the RESD design amount to \SI{28.0}{\million\euro\per\year} and are comprised of \SI{19.2}{\million\euro\per\year} of capital expenses and \SI{8.8}{\million\euro\per\year} of operational expenses.
The average cost of electricity generation is \SI{105.6}{\euro\per\mega\watt\hour}, compared to the actual average cost of electricity generation in the entire Canary Islands of \SI{161.5}{\euro\per\mega\watt\hour} in \si{2021} \citep{gobiernodecanariasAnuarioEnergeticoCanarias2023}.
Note that the current energy system has an overcapacity of conventional diesel generators due to some of the generation capacity being held back as a backup.
In the RESD approach, we do not take into account possible component failures.

In the robust design, the majority of electricity is produced by wind turbines with \SI{50.8}{\percent}, followed by solar with \SI{41.2}{\percent}, and diesel with \SI{8.0}{\percent}.
Renewable resources are preferred over diesel since they have lower variable operational costs.
The high variable operational costs of diesel engines are primarily driven by the price of fuel and the costs of dispatch emission rights, which are described in Section 2.3 of the supplementary material.
Compared to the current design, the RESD design exhibits a much higher share of renewable generation, with the renewable energy penetration reaching \SI{92.0}{\percent}.
In comparison, only \SI{10.3}{\percent} of gross electricity production in La Palma was supplied by renewable sources in 2021 \citep{gobiernodecanariasAnuarioEnergeticoCanarias2023}.

\subsection{Performance and Accuracy}
We run the La Palma case study for different numbers of time steps and principal components (PCs). Varying the latter parameter allows us to determine how the dimensionality reduction influences the computational performance and the accuracy of the obtained results.
We always use five random seeds for the Gurobi optimizer \citep{gurobi} and average the obtained solution times since solution times can vary strongly with random seeds depending on the nodes explored by the branch and bound algorithm.
We measure the solution time by taking the `CPU time' reported by libDIPS  \citep{zinglerLibDIPSDiscretizationBasedSemiInfinite2023}, which is the sum of the wall-clock times reported by Gurobi \citep{gurobi} for all sub-problems.
Note that we did not account for the time it takes to calculate the convex hull in the performance evaluations.
As mentioned in Section~\ref{sec:data_processing}, with larger latent space dimensionality, the time to compute the convex hull can become significant compared to the solution time.
Hence, we limited the maximum number of principal components in our numerical experiments to \si{9}.

\begin{figure}[H]
    \centering
    \includegraphics[width=222pt]{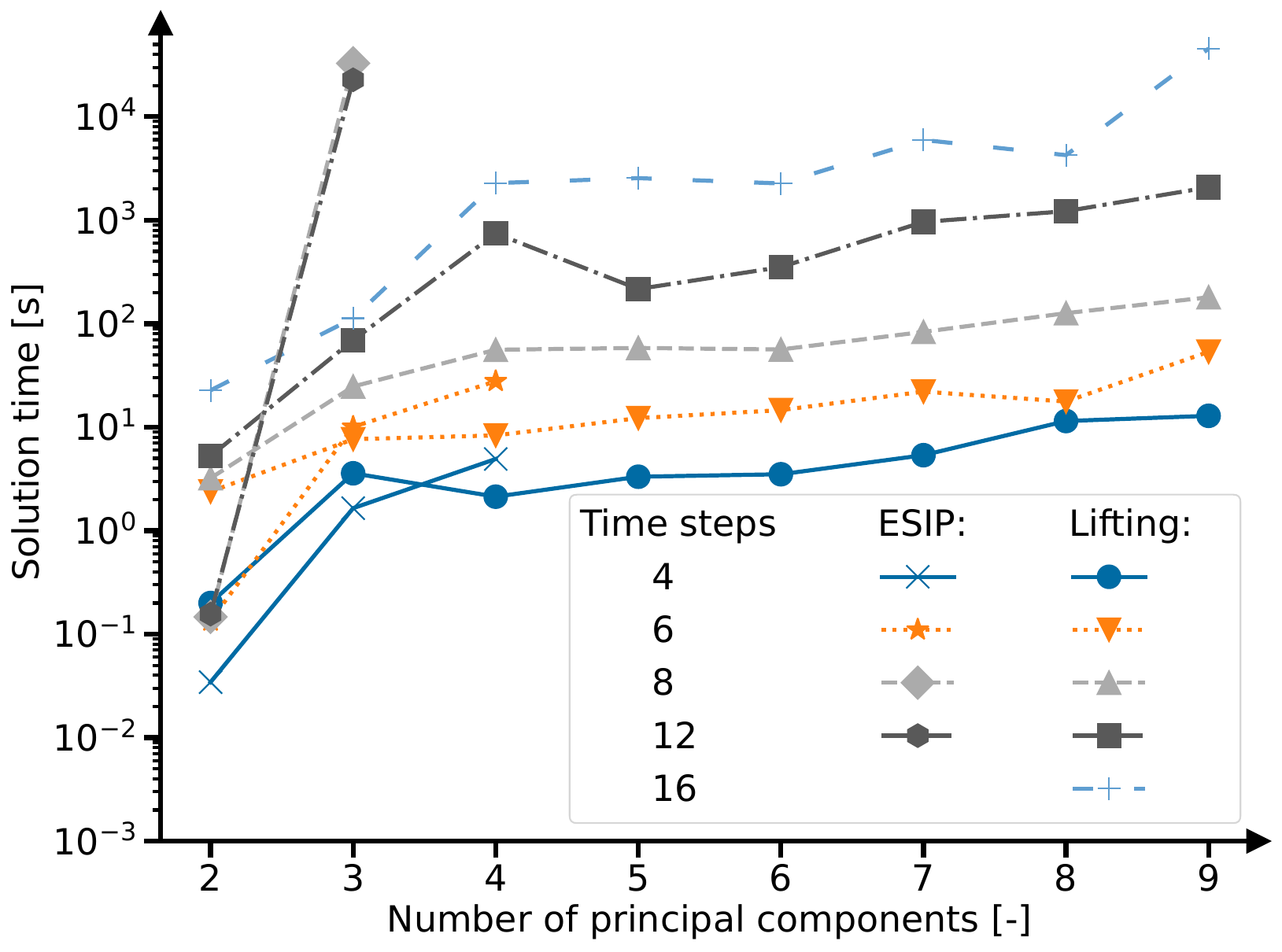}
    \caption{Average solution times in logarithmic scale for the ESIP and the lifting approach. The average solution time in seconds is plotted against the number of principal components (PCs) for varying numbers of time steps. The solution time increases both with an increasing number of PCs and time steps. The lifting approach scales better with an increasing number of PCs than the ESIP approach.}
    \label{fig:performance}
\end{figure}

Figure \ref{fig:performance} shows the average solution time for the ESIP and the lifting approach for different numbers of PCs and time steps.

As expected, the solution time increases with an increasing number of time steps and PCs.
The ESIP approach shows a very strong increase in the solution time with an increasing number of PCs.
In fact, we could not obtain results for four or more components with the ESIP approach since the run times became too long.
The lifting approach outperforms the ESIP approach if more than two PCs are used.
There is a significant increase in solution time between two and three components.
Afterward, solution time increases more slowly with an increasing number of components.
Furthermore, higher temporal resolutions strongly increase solution times.

\begin{figure}[H]
    \centering
    \includegraphics[width=222pt]{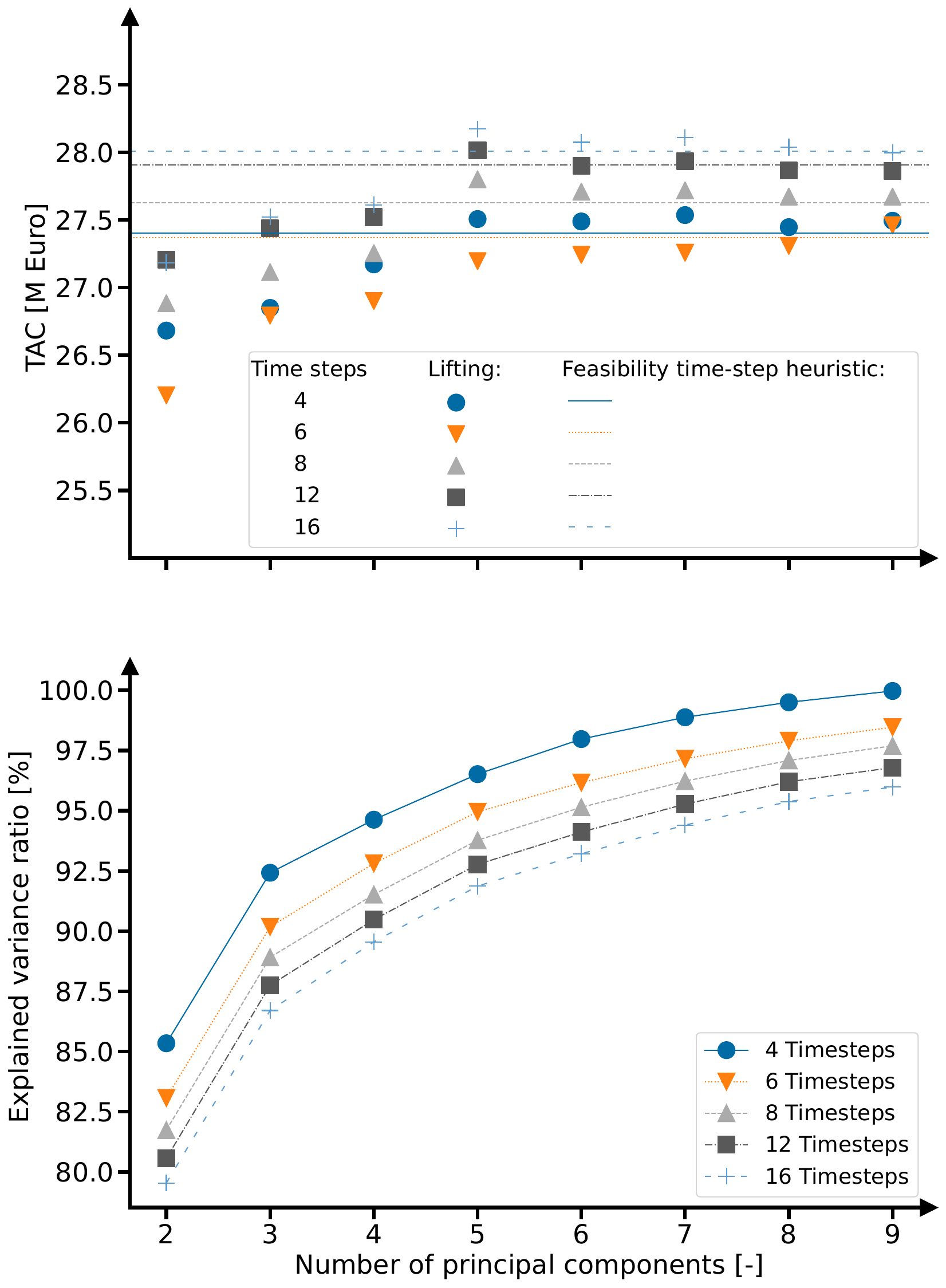}
    \caption{Optimal total annualized cost (TAC) of energy system designs obtained using the lifting approach against the number of principal components (PCs) for varying number of time steps (top). Dashed lines show results obtained with the feasibility time-step heuristic \citep{bahlTimeseriesAggregationSynthesis2016, teichgraeberExtremeEventsTime2020} discussed in Section \ref{special}. For a low number of PCs, a significant approximation error can be seen. If the number of PCs is chosen sufficiently large, i.e., greater than or equal to $5$, the TAC of full-dimensional problems is approximated closely by the RESD approach. Explained variance ratio plotted against the number of PCs for a varying number of time steps (bottom). A small number of PCs can explain a majority of the variance in the historical data.}
    \label{fig:accuracy}
\end{figure}

The top of Figure \ref{fig:accuracy} shows the optimal TACs of the energy system designs obtained using the lifting approach plotted against the number of PCs for a varying number of time steps.
As a reference, the results obtained using the feasibility time-step heuristic \citep{bahlTimeseriesAggregationSynthesis2016, teichgraeberExtremeEventsTime2020} described in Section \ref{special} are shown.
The bottom shows the explained variance ratio.

If the number of PCs is small, the RESD designs have lower TACs than the designs resulting from the heuristic.
This indicates that, for a small number of PCs, the truncation error from the dimensionality reduction leads to designs that are not robust.
With an increasing number of PCs, the TAC of the RESD designs quickly approach those of the designs resulting from the heuristic.
In some cases, the TACs of the RESD designs exceed those of the designs obtained by the heuristic, indicating that the truncation error can also lead to overly robust designs.
Note that relatively few PCs are needed to account for a majority of the variance in the historical data.
With \si{5} PCs, more than \SI{90}{\percent} of the variance are explained for all time resolutions.
We observe that the TAC of the full-dimensional problem are approximated closely for all time resolutions if the number of PCs is chosen such that the explained variance is above \SI{95}{\percent}.

\begin{figure}[H]
    \centering
    \includegraphics[width=222pt]{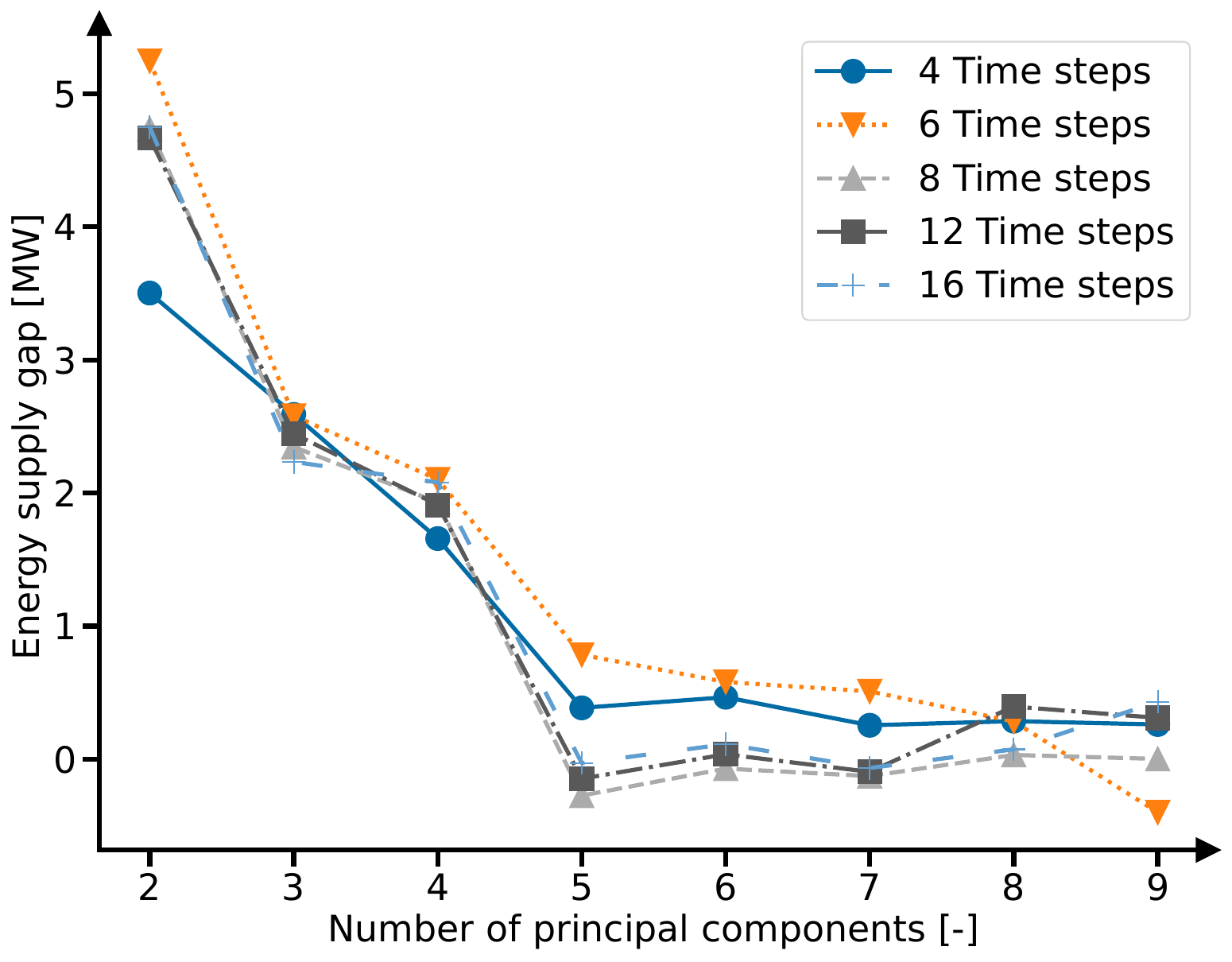}
    \caption{The maximum energy supply gap for RESD designs with different time resolutions obtained using the lifting approach, calculated by solving the operational problem for all historical data points, is plotted against the number of principal components (PCs). The supply gap is large for a small number of PCs, showing that these designs are not robust. The energy supply gap decreases quickly with an increase in the number of PCs and levels off close to \SI{0}{\mega\watt} at \si{5} PCs, indicating robust designs.}
    \label{fig:supply_gap}
\end{figure}

Figure \ref{fig:supply_gap} shows the energy supply gap over the number of PCs for different time resolutions.
The energy supply gap is calculated by solving the operational problem for each historical data point and selecting the maximum constraint violation.

For a small number of PCs, the energy supply gap is large, i.e., in the order of multiple megawatts, showing that these designs are not robust.
This finding highlights the tradeoff between a sufficiently accurate approximation of the time-series data and the computational performance.
Importantly, the energy supply gap decreases quickly with an increase in the number of PCs and levels off around \si{5} PCs.
In fact, for \si{48} dimensional data, i.e., \si{16} time steps and \si{3} quantities (solar, wind, demand), \si{6} PCs are sufficient to achieve an appropriate approximation of the full-dimensional problem and hence a robust design.
In some cases, the energy supply gap even becomes negative, indicating an overly robust design.
This is in accordance with the observation that the TACs for these numbers of PCs are higher than the respective TACs for the designs obtained by the heuristic.

Choosing the right dimensionality of the latent space is difficult in the case of a nonconvex problem.
However, we note that the energy supply gap approaches zero for designs that were obtained with a latent space dimensionality that covers more than \SI{95}{\percent}
of the variance in the historical data.
Hence, the explained variance ratio may be used as an indicator to decide on an appropiate dimensionality of the latent space.

\section{Conclusion}\label{Conclusion}
Uncertainties introduced by VRES should be incorporated into the design of energy systems.
Previously \emph{a priori} heuristics, e.g., statistical selection of extreme periods \citep{dominguez-munozSelectionTypicalDemand2011}, or optimization-based heuristics, e.g., the feasibility time-step heuristic \citep{bahlTimeseriesAggregationSynthesis2016, teichgraeberExtremeEventsTime2020}, have been proposed.
However, \emph{a priori} methods cannot identify extreme periods specific to the energy system because they are agnostic to its design.
Moreover, optimization-based heuristics may fail if the design problem exhibits a nonconvex operational problem, e.g., a MILP resulting from the need to model minimal part loads.

We present the robust energy system design (RESD) approach, a rigorous but computationally intense approach to identify extreme scenarios during the design optimization.
Specifically, we introduce a semi-infinite existence constraint to the energy system design problem to ensure  robustness with respect to uncertainties, resulting in a tri-level optimization problem.

Our RESD approach for the identification of worst-case scenarios during optimization can be seen as a generalization of the finite optimization problem corresponding to the feasibility time-step heuristic \citep{bahlTimeseriesAggregationSynthesis2016, teichgraeberExtremeEventsTime2020}, which determines the constraint violation for all historical scenarios and iteratively adds scenarios with infeasible time steps as extreme scenarios.
However, as we have discussed, it is only sufficient to check operability with the historical data, or more specifically, all the scenarios that define the boundary of the convex hull of the historical data, in the case of linear medial- and lower-level problems or jointly convex medial- and lower-level problems.
The computational performance of the iterative feasibility time-step heuristic could potentially be improved by using a preprocessing step to identify and discard some or all data points in the interior. In contrast, our RESD approach is dedicated to the more general nonconvex case. 
The MILP example in Section \ref{MILP Example} demonstrates how components with minimal part load in combination with a no-curtailment constraint can lead to worst-case scenarios that lie between historical data and not on the vertices of the uncertainty set.
Future work should examine to which extent the worst-case scenarios for realistic multi-energy systems, which allow for some curtailment and/or contain storage components, also do not lie on the vertices of the uncertainty set.

Different uncertainty-bounding methods can be used in the RESD approach.
In the present contribution, we consider the convex hull of the historical data as the possible uncertainty set, taking advantage of the fact that the describing equations are linear.
Still, the computational effort for solving the RESD problem increases quickly with an increasing number of time steps.

Computational tractability can be improved by leveraging dimensionality reduction techniques such as principal component analysis (PCA).
While we observe that robust designs can be obtained even in the presence of rather aggressive dimensionality reduction, the use of PCA inherently introduces an approximation error that can lead to overly or underly robust designs. A suitable dimensionality of the latent space must be found empirically, and in the nonconvex case, a rigorous criterion for validating the robustness of an obtained design is missing.
However, the explained variance ratio can be used as an indicator in guiding the selection of an appropriate latent space dimensionality.
Specifically, we found designs obtained with a latent space dimensionality covering more than \SI{95}{\percent} of the variance of the historical data to yield TACs that are close to the TAC of the full-dimensional design.

Although the lifting approach can improve performance if the considered design problem has an embedded convex operational problem, the applicability of the RESD approach is currently limited to small problems due to the involved computational intensity.
Furthermore, we found that including the complementarity constraints improves the performance of the lifting approach.
This is an unexpected result, and it should be investigated if this holds generally.
Wider applicability of our approach could be enabled by faster methods for solving hierarchical programs.
For example, the approach by \cite{seidelAdaptiveDiscretizationMethod2022} that improves the rate of convergence of the adaptive discretization approach could reduce the long solution times of the embedded MaxMin problem.
Furthermore, the integration of bounding approaches that can handle nonconvexities in the historical data, e.g., holes, as well as the use of nonlinear dimensionality reduction techniques, e.g., autoencoders \citep{kramerNonlinearPrincipalComponent1991}, should be investigated.

\section*{Declarations}
\subsection*{Ethics Approval and Consent to Participate}
Not applicable.

\subsection*{Consent for Publication}
Not applicable.

\subsection*{Funding}
This work was performed as part of the Helmholtz School for Data Science in Life, Earth and Energy (HDS-LEE) and received funding from the Helmholtz Association of German Research Centers.

The authors gratefully acknowledge the financial support of the Kopernikus project SynErgie by the German Federal Ministry of Education and Research (BMBF) and the project supervision by the project management organization Projektträger Jülich (PtJ).

\subsection*{Availability of Data and Materials}
The cost parameters for the La Palma case study are collected in Table \si{1} in the supplementary material, and a detailed description of how they were obtained can be found in Section $2.3$ of the supplementary material.
The solar irradiance and wind speed data used in the La Palma case study are available at \url{https://re.jrc.ec.europa.eu/pvg_tools/en/}.
The demand data for La Palma island is available at \url{https://demanda.ree.es/visiona/canarias/la_palma5m/tablas}.
The wind turbine power curve data for an Enercon E-82 E2 turbine was used.
While the data is no longer available on the official manufacturer site, a copy of the data sheet is accessible at \url{https://catalystresearch.wordpress.com/wp-content/uploads/2013/11/enercon_pu_en.pdf}.

\subsection*{Competing Interests}
The authors have no relevant financial or non-financial interests to disclose.

\section*{Authors' Contributions}
Conceptualization: M.W., E.C., A.M., M.D.; Methodology: M.W.; Software: M.W.; Formal analysis and investigation: M.W.; Visualization: M.W.; Writing - original draft preparation: M.W.; Writing - review and editing: E.C., A.M., M.D.; Funding acquisition: A.M., M.D.; Supervision: A.M., M.D.

\section*{Acknowledgements}
During the preparation of this work, M.W. used Grammarly to correct grammar and spelling and improve the writing style. After using this tool, all authors reviewed and edited the content as needed and take full responsibility for the content of the publication.

\section*{Nomenclature}
Throughout the manuscript, scalar-valued quantities are denoted in regular font, e.g., $x$, vector-valued quantities are denoted in bold font, e.g., $\mathbf{x}$, and set-valued quantities are denoted in calligraphic font, e.g., $\mathcal{X}$.
\subsection*{Abbreviations}
\begin{table}[H]
\begin{tabular}{ll}
    BLLP & bi-level linear program \\
    ESIP & existence-constrained semi-infinite program \\
    EGSIP & existence-constrained generalized semi-infinite program \\
    GSIP & generalized semi-infinite program \\
    MILP & mixed-integer linear program \\
    MW & megawatt \\
    NLP & nonlinear program \\
    RESD & robust energy system design \\
    PC & principal component \\
    PCA & principal component analysis \\
    PV & photo voltaic \\
    SIP & semi-infinite program \\
    TAC & total annualized cost \\
    VRES & variable renewable energy sources \\
\end{tabular}
\end{table}

\subsection*{Greek Symbols}
\begin{table}[H]
\begin{tabular}{ll}
    $\lambda$ & Lagrange multiplier for equality constraints \\
    $\mu$ & Lagrange multiplier for inequality constraints \\
    $\phi$ & auxiliary variable \\
\end{tabular}
\end{table}

\subsection*{Latin Symbols}
\begin{table}[H]
\begin{tabular}{ll}
    $b$ & binary variable for on/off decision \\
    $e$ & energy supply gap \\
    $e_{epi}$ & auxiliary variable for energy supply gap \\
    $g$ & inequality constraint \\
    $\mathbf{g}$ & vector of inequality constraints \\
    $h$ & equality constraint \\
    $\mathbf{h}$ & vector of equality constraints \\
    $\mathcal{L}$ & Lagrangian function \\
    $n$ & number\\
    $\mathcal{S}$ & set of representative scenarios \\
    $\mathcal{T}$ & set of time steps \\
    $\mathbf{w}$ & vector of variables for single-level NLP \\
    $\mathcal{W}$ & feasible set of variables for single-level NLP \\
    $x$ & design variable \\
    $\mathbf{x}$ & vector of design variables \\
    $\mathcal{X}$ & feasible set of design variables \\
    $y$ & uncertain variable \\
    $\mathbf{y}$ & vector of uncertain variables \\
    $\tilde{\mathbf{y}}$ & vector of dependent uncertain variables \\
    $\bar{\mathbf{y}}$ & vector of independent uncertain variables \\
    $\mathcal{Y}$ & feasible set of uncertain variables\\
    $z$ & operational variable \\
    $\mathbf{z}$ & vector of operational variables \\
    $\mathcal{Z}$ & feasible set of operational variables \\
\end{tabular}
\end{table}

\subsection*{Subscripts}
\begin{table}[H]
\begin{tabular}{ll}
    $d$ & data points\\
    $e$ & energy supply gap \\
    $en$ & energy system \\
    $eq$ & equations \\
    $i$ & index \\
    $j$ & index \\
    $epi$ & epigraph reformulation \\
    $o$ & operational costs \\
    $p$ & principal components \\
    $ref$ & reformulation \\
    $s$ & representative scenarios \\
    $t$ & time step \\
    $x$ & design variables \\
    $y$ & uncertain variables \\
    $z$ & operational variables \\
\end{tabular}
\end{table}
 
\bibliographystyle{apalike}
  \renewcommand{\refname}{References}  
  \bibliography{bibliography.bib}
\end{document}


\thispagestyle{firststyle}

  \begin{center}
    \begin{large}
      \textbf{\mytitle}
    \end{large} \\
    \myauthor
  \end{center}

  \vspace{0.5cm}

  \begin{footnotesize}
    \affil
  \end{footnotesize}

  \vspace{0.5cm}

\section{Illustrative Nonconvex Example: ESIP Subproblems}\label{MILP_Subproblems}
This section complements Section \si{3} of the main manuscript and provides the subproblems involved in the solution of the mixed-integer linear program (MILP) existence-constrained semi-infinite programming (ESIP) problem.
First, we have the lower bounding problem:
\begin{alignat*}{4}\tag{MILP LBP}\label{MILP_lbp}
    &\underset{x_1, x_2, z_{1, u, k}, b_{u, k}}{\min} && 2x_1 + x_2\\
    &\text{s.t.} && y_{k, disc} - z_{1, u, k} - x_2 b_{u, k} \le 0 \quad \forall k \in \mathcal{K}\\
    &&&  0.2 x_2 b_{u, k} + z_{1, u, k} - y_{k, disc} \le 0 \quad \forall k \in \mathcal{K}\\
    &&& z_{1, u, k} - x_1 \le 0 \quad \forall k \in \mathcal{K}\\
    &&& x_1, x_2, z_{1, u, k} \in [0, 100]\\
    &&& b_{u, k} \in \{0, 1\}
\end{alignat*}
The lower bounding problem \eqref{MILP_lbp} determines a lower bound on the optimal objective value and fixes the capacities for the subsequent solution of the MAXMIN problem that determines whether the identified design is feasible.
The discretization set $\mathcal{K}$, populated by the medial-level problem, comprises the identified extreme scenarios up to this point and ensures that the identified design is robust.
$x_1$ and $x_2$ are abstract component capacities.
$z_{1, u, k}$ is the operational variable of the lower bounding problem.
There are $|\mathcal{K}|$ sets of operational variables, one for each discretization point.
$b_{u, k}$ are the lower bounding binary variables used to model the minimal part load of Component 2.
The discretization parameter $y_{k, disc}$ corresponds to the optimal solution in iteration $k$ to the medial-level problem:
\begin{alignat*}{4}\tag{MILP MLP}\label{MILP_mlp}
    &\underset{y, mlpobj, b_m}{\min}&& -mlpobj\\
    &\text{s.t.}&& \ mlpobj - y + z_{1, l, disc} + b_{l, disc} x_2 - 1000b_m \le 0 \quad  \forall l \in \mathcal{N}\\
    &&&  mlpobj - b_{l, disc} 0.2 x_2 + y - z_{1, l, disc}- 1000(1 - b_m) \le 0 \quad  \forall l \in \mathcal{N}\\
    &&& y \in [0, 100]\\
    &&& mlpobj \in [-1000, 1000]\\
    &&& b_{m} \in \{0, 1\}
\end{alignat*}
The medial-level problem \eqref{MILP_mlp} determines worst-case demand realizations for given component capacities $x_1$ and $x_2$ for a given discretization set $\mathcal{N}$ populated by the lower-level problem.
$y$ is the uncertain demand, and $mlpobj$ is a dummy variable introduced by the epigraph reformulation of the medial-level problem \citep{djelassiDiskretisierungsbasierteAlgorithmenFur2020}.
$b_m$ is a binary variable introduced to reformulate the maximum introduced by moving the lower-level constraints depending on $y$ into the objective.
Discretization parameters $z_{1, l, disc}$ and $b_{l, disc}$ correspond to the optimal values of the operational variables and are determined by the lower-level problem:
\begin{alignat*}{4}\tag{MILP LLP}\label{MILP_llp}
    &\underset{z_{1, ll}, b_{ll}, llpobj}{\min} && llpobj\\
    &\text{s.t.} && y - z_{1, ll} - b_{ll} x_2 - llpobj \le 0\\
    &&& 0.2 b_{ll} x_2 - y + z_{1, ll} - llpobj \le 0\\
    &&& z_{1, ll} - x_1 \le 0\\
    &&& z_{1, ll} \in [0, 100]\\
    &&& llpobj \in [-1000, 1000]\\
    &&& b_{ll} \in \{0, 1\}
\end{alignat*}
The lower-level problem \eqref{MILP_llp} describes the operational problem for fixed capacities $x_1$ and $x_2$ and fixed demand $y$.
The lower-level problem minimizes the constraint violations of the operational problem and determines the discretization points for the medial-level problem \eqref{MILP_mlp}.
Operational variables $z_{1, ll}$ and $b_{ll}$ are the lower-level operational variables.

\section{La Palma Energy System}\label{sec:CaseStudy}
This section provides further details on the data processing as well as the complete model formulations for the La Palma energy system introduced in Section \si{4} of the main manuscript.

\subsection{Historical Time-Series Data}
For the renewable energy resources, we use data for wind speed at \SI{10}{\meter} height and global irradiance obtained from the Photovoltaic Geographical Information System (PVGIS) by the European Union \citep{europeancommissionJRCPhotovoltaicGeographical2022, huldNewSolarRadiation2012}.
We chose the northwest of La Palma (\SI{28.780}{\degree} latitude and \SI{-17.980}{\degree} longitude) as the location for the historical data after a visual inspection using PVGIS \citep{europeancommissionJRCPhotovoltaicGeographical2022, huldNewSolarRadiation2012} showed both annual average solar irradiance and average wind speed to be highest in the northwest \citep{europeancommissionJRCPhotovoltaicGeographical2022, davisGlobalWindAtlas2023}.
PVGIS uses SARAH2.1 \citep{pfeifrothSurfaceRadiationData2019} and ERA5 \citep{hersbachERA5GlobalReanalysis2020} as underlying data sets for solar irradiance and wind speed, respectively.
The solar panels are assumed to have an optimal orientation, which PVGIS determines for the chosen location.

The PV unit and wind turbine power outputs are modeled using utilization factors $f_{c, s, t}$, which are bounded by capacity factors $c_{cap, solar, s, t}$ calculated during a preprocessing step.
For PV power, the solar capacity factor is calculated from the global irradiance according to \cite{sassModelCompendiumData2020} as
\begin{equation*}
     c_{cap, solar, s, t} = \min\{ \frac{I_{s, t} \eta_{solar}}{P_{solar, nom}}, 1\},
\end{equation*}
where $I_{s, t}$ is the global irradiance in \SI{}{\kilo\watt\per\square\meter}, $\eta_{solar} = 0.19$ is the assumed efficiency of the PV units, and $P_{solar, nom} = \SI{0.171}{\kilo\watt\per\square\meter}$ is the nominal capacity.\\
The wind capacity factor is calculated from the wind speed at hub height using the power curve of an Enercon E82/2350 wind turbine \citep{Enercon_2012} as
\begin{equation*}
    c_{cap, wind, s, t} = \frac{P_{wind, s, t}}{P_{wind, nom}},
\end{equation*}
where $P_{wind, nom} = \SI{2350}{\kilo\watt}$ is the nominal capacity of a wind turbine.
The power output $P_{wind, s, t}$ for a given wind speed at hub height $v_{wind, hub, t}$ in \si{\meter\per\second} is calculated by linear interpolation of the power output from the power curve data, which is given for wind speeds in \SI{1}{\meter\per\second} increments.
We obtain wind speed at hub height from wind speed at \SI{10}{\meter} by using the logarithmic wind profile approximation \citep{okeBoundaryLayerClimates2006}
\begin{equation*}
      v_{wind, hub} = v_{wind, 10m} \frac{\ln \left( h_{hub} / z_0 \right)}{\ln \left( \SI{10}{m} / z_0 \right)},
\end{equation*}
with a hub height $h_{hub}$ of \SI{85}{\meter} and a roughness length $z_0$ of \SI{0.3}{\meter}, determined by visual inspection of the roughness length map in PVGIS \citep{europeancommissionJRCPhotovoltaicGeographical2022, huldNewSolarRadiation2012} for the north-west of La Palma \citep{davisGlobalWindAtlas2023}.
The cutout wind speed for the turbine is \SI{25}{\meter\per\second}.
At higher wind speeds, the turbine is shut down, and no power is produced.

\begin{figure}[H]
    \centering
    \includegraphics[width=255pt]{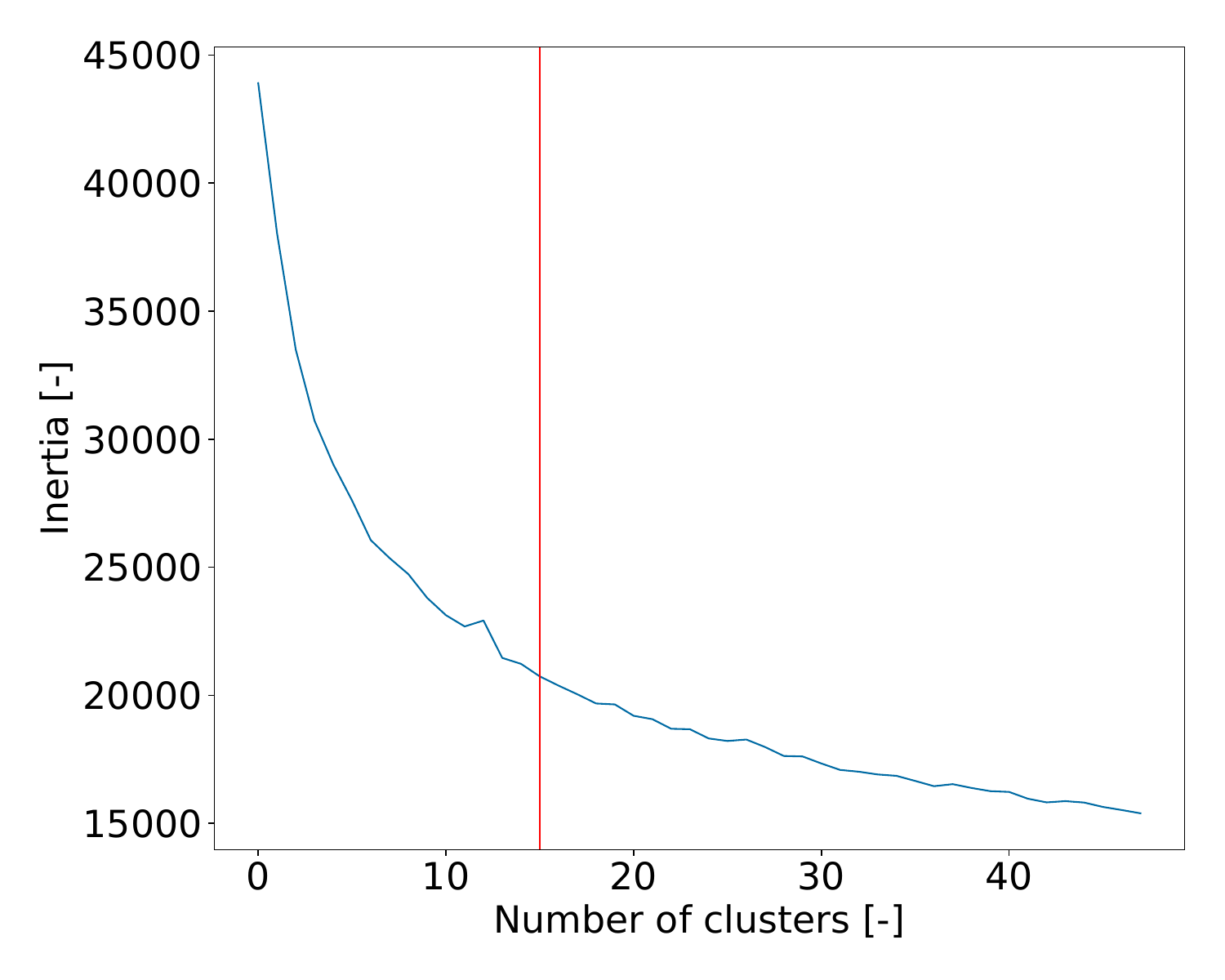}
    \caption{k-means clustering inertia plotted against the number of clusters (elbow method \citep{thorndikeWhoBelongsFamily1953}). Clusters are determined using the scikit-learn package \citep{scikit-learn}, default settings are used, and the random state is fixed to $42$ to ensure reproducibility. The red line indicates the number of clusters chosen in this study, which is \si{15}.}
    \label{fig:elbow}
\end{figure}

\begin{figure}[H]
    \centering
    \includegraphics[width=255pt]{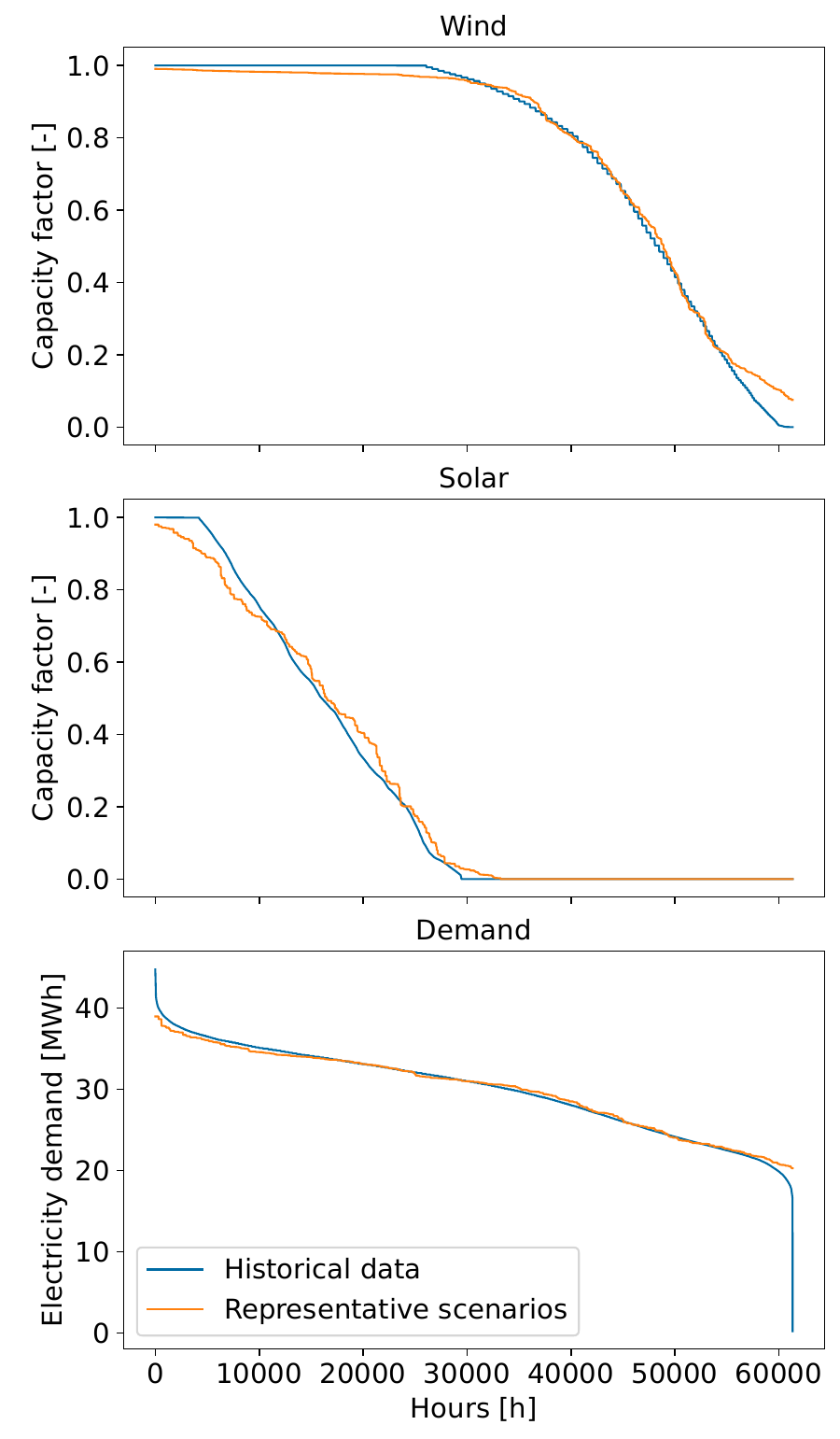}
    \caption{Load curves \citep{mastersRenewableEfficientElectric2004} for the considered time-series data. Historical data and data of the 15 selected representative days. To compute the load curve for the representative days, we consider the full historical time series but replace each historical day by the corresponding representative day. A good agreement between the load curves of the historical data and the representative scenarios can be seen, indicating that the representative scenarios can approximate the total yearly load well.}
    \label{fig:load_curve}
\end{figure}

Demand data is aggregated from the Spanish electricity distribution system operator \cite{redelectricadeespanaPalmaElectricityDemand2024}.
Solar, wind, and demand data covers the time span from 01.01.2013 until 31.12.2019.
We follow the processing introduced in Section \si{2} of the main manuscript to obtain the representative days and the uncertainty bounds.
First, z-normalization is applied to each individual quantity, i.e., solar capacity factors, wind capacity factors, and electricity demand \citep{teichgraeberClusteringMethodsFind2019}.
Normalization is performed for the full time-series data, i.e., for each quantity, all values in the historical data are considered at once, as opposed to normalizing each representative period or each time step individually.
The normalized time series data for each attribute are then concatenated, and representative days are selected by k-means clustering \citep{macqueen1967some} using the scikit-learn package \citep{scikit-learn}.
Default settings are used and the random state is fixed to $42$ to ensure reproducibility.
To accurately approximate the operational costs, an appropriate number of clusters and, thus, representative scenarios must be selected.
The elbow method \citep{thorndikeWhoBelongsFamily1953}, shown in Figure \ref{fig:elbow}, does not clearly indicate a suitable number of representative days.
Furthermore, the elbow method has been criticized for its strong dependence on the scaling of the plot axis and spurious results in case of uncorrelated data \citep{schubertStopUsingElbow2023}.
Pragmatically, we choose 15 representative days to achieve a reasonable balance between the computational effort and the accuracy of the load curve approximation (see Figure \ref{fig:load_curve}).


To determine the uncertainty bounds, we reduce the dimensionality of the concatenated time-series data using principal component analysis (PCA) \citep{pearsonLinesPlanesClosest1901} implemented in scikit-learn \citep{scikit-learn}.
We use 'full' for the svd\_solver parameter to calculate the exact PCA decomposition and default values for all remaining parameters.
Then, we determine the convex hull, using the SciPy package \citep{virtanenSciPyFundamentalAlgorithms2020} on the latent space historical data.

\subsection{Model Formulation}\label{sec:Problem Formulation}
For ease of notation, we define the set of generation components $\mathcal{C} = \{solar, wind, diesel\}$.
Similar to \cite{sassModelCompendiumData2020}, we use total annualized cost (TAC) as our objective, i.e., 
\begin{align*}
    \underset{\mathbf{x} \in \mathcal{X}, \mathbf{z}_s \in \mathcal{Z}_s}{\min} \sum_{c\in \mathcal{C} \cup \{battery\}} P_{c, peak} \left(\frac{c_{inv, c}}{f_{AN}} + c_{fix, c}\right) \\ + 365 \sum_{s \in \mathcal{S}} \omega_s \sum_{c\in \mathcal{C}}{ \sum_{t \in \mathcal{T}} \Delta_t P_{c, peak} f_{c, s, t} c_{var, c}}, \label{eqn:obj}
\end{align*}
with design decisions $\mathbf{x}$ and operational variables $\mathbf{z}_s$.
The objective function consists of investment costs for each generation component and the battery $c$, calculated by multiplying its installed capacity $P_{c, peak}$ with the specific investment costs  $c_{inv, c}$.
Similarly, fixed operational costs $c_{fix, c}$ are considered for each year of operation.
Additionally, investment costs are divided by the present value annuity factor
\begin{equation*}
    f_{AN} = \frac{(1 + i)^T - 1}{i (1 + i)^T},
\end{equation*}
which is determined with an assumed interest rate $i$ of \SI{8}{\%} and a time horizon $T$ of \si{25} years.
Finally, variable operational costs $c_{var, c}$ are multiplied by the provided amount of electricity by each generation component in each time step $t$ for every representative scenario $s$.
To this end, the length of each time step $\Delta_t$ is multiplied by the installed capacity and the capacity factor $f_{c, s, t}$.
The cost for each representative scenario are then weighted by $\omega_s$, the number of days in each cluster $n_{days, s}$ divided by the number of total days $n_{days}$, see e.g., \cite{nahmmacherCarpeDiemNovel2016}, i.e., 
\begin{equation*}
    \omega_s = \frac{n_{days, s}}{n_{days}} .
\end{equation*}
Yearly variable operational costs are calculated by scaling the costs of the representative days to a year.
\begin{table}
    \begin{center}
        \begin{tabular}{ c|ccc } 
         \hline
         Component $c$ & $c_{inv, c}$ in \si{\euro_{2022}\per\kilo\watt} & $c_{fix, c}$ in \si{\euro_{2022}\per\kilo\watt} & $c_{var, c}$ in \si{\euro_{2022}\per\kilo\watt\hour}\\
         \hline
         Solar & 883.3  & 17.9 & 0 \\ 
         Wind & 2283.7  & 26.9  & 0.011 \\
         Diesel & 2391.8 & 0 & 0.242 \\
         Battery & 1550 & 31 & 0 \\
         \hline
        \end{tabular}
    \caption{Cost parameters of the considered components.\label{tab:cost}}
    \end{center}
\end{table}
Table \ref{tab:cost} summarizes all cost parameters used in this study.
We adjusted all costs to a 2022 base using the EU domestic producer price index for the 19-country Euro area \citep{eurostatDomesticProducerPrices2024}.
Details on how the cost parameters are obtained and adjusted for inflation are given in Section \ref{Parameters}.

We introduce a constraint to ensure that the power demand is satisfied for every time step in every scenario: 
\begin{equation*}
    P_{demand, s, t} + P_{battery, in, s, t} - P_{battery, out, s, t} - \sum_{c\in \mathcal{C}} P_{c, peak} f_{c, s, t} \le 0 \ \forall s \in \mathcal{S} \ \forall t \in \mathcal{T}
\end{equation*}
$P_{demand, s, t}$ is the electricity demand at time step $t$ of a representative day $s$.
Equation $f_{c, s,t} - c_{cap, c, s, t} \le 0$ describes that capacity factors $f_{c, s, t}$ have to be less than the capacity factor parameters for the representative scenarios $c_{cap, c, s, t}$.
For the diesel generator the capacity factor parameters are \si{1}, i.e., $c_{cap, diesel, s, t} = 1$.

The battery system is modeled following \cite{sassModelCompendiumData2020}.
Separate variables $P_{battery, in, s, t}$ and $P_{battery, out, s, t}$ are introduced for charging and discharging to account for charging and discharging losses.
The relation between peak power output $P_{battery, peak}$ and peak capacity $E_{battery, peak}$ is assumed to be \SI{4}{\kilo\watt\hour\per\kilo\watt}, i.e., a constraint $E_{battery, peak} = 4 P_{battery, peak}$ is added.
The change in the stored electricity $E$ is described by
\begin{equation*}
    \frac{dE}{dt} = \frac{1}{3600} (\eta_{in} P_{battery, in} + \frac{1}{\eta_{out}} P_{battery, out}),
\end{equation*}
with charging efficiency $\eta_{in} = 0.92$ and discharging efficiency $\eta_{out} = 0.926$ based on \cite{baumgartnerDesignLowcarbonUtility2019}.
We discretize the differential equation using the forward Euler scheme \citep{eulerFoundationsDifferentialCalculus2000}, i.e., we obtain
\begin{align*}
    \dot E_{battery, s, t} - \frac{1}{3600} (\eta_{in} P_{battery, in, s, t} + \frac{1}{\eta_{out}} P_{battery, out, s, t}) = 0 \ \forall s \in \mathcal{S} \ \forall t \in \mathcal{T},\\
    \dot E_{battery, s, t} \Delta t - (E_{battery, s, t} - E_{battery, s, t-1}) = 0 \ \forall s \in \mathcal{S} \ \forall t \in \mathcal{T},
\end{align*}
where $\Delta t$ is the time step length, $E_{bat, s, t}$ is the state of charge, and $\dot E_{bat, s, t}$ the rate of charge.
The initial state of charge is assumed to be \SI{50}{\%} of the nominal capacity.
Additionally, there is a cyclical constraint requiring that the battery is recharged to the initial level at the end of the representative period.

We model the diesel generator as a dispatchable electricity source, with power output calculated from installed capacity $P_{diesel, peak}$ and variable utilization factor $f_{diesel, s, t}$.
We assume that ramping constraints for the diesel generator can be neglected for the investigated temporal resolution.
Furthermore, we assume that no minimal part load is needed since the installed diesel generators consist of multiple engines whose individual part loads are small compared to the overall load.

An embedded optimization problem is added that guarantees the feasibility of the design for the worst-case uncertainty realization:
\begin{equation*}
     \underset{\mathbf{y}\in \mathcal{Y}(\mathbf{x})}{\max}\underset{\mathbf{z}\in \mathcal{Z}_{epi}(\mathbf{x}, \mathbf{y})}{\min} \underset{t\in \mathcal{T}}{\max} \ P_{demand, wc, t} + P_{battery, in, wc, t} - P_{battery, out, wc, t} - \sum_{c\in \mathcal{C}} P_{c, peak} f_{c, wc, t}  \le 0
\end{equation*}
Its objective is to maximize the electricity supply gap by searching for the worst-case renewable capacity factors $f_{wind, wc, t}$, $f_{solar, wc, t}$, and demand realizations $P_{demand, wc, t}$.
To compensate for the worst-case realization, the battery operation, i.e., $P_{battery, in, wc, t}$ and $P_{battery, out, wc, t}$, as well as the diesel engine utilization factor $f_{diesel, wc, t} $ can be adjusted.
Notice that the utilization factor for the diesel engine $f_{diesel, wc, t}$ can be chosen freely within the interval $[0, 1]$ since it does not appear in any of the lower-level constraints.
It has a negative sign in the objective function of the embedded minimization problem, so an optimal solution exists where it is at its upper bound $1$. Thus, we fix $f_{diesel, wc, t} = 1$ to reduce the computational complexity.

We use PCA to reduce the dimensionality of the optimization problem to facilitate computational tractability.
To this end, we calculate the worst-case demand $P_{demand, wc, t}$, as well as the worst-case wind and solar capacity factors $f_{wind, wc, t}$ and $f_{solar, wc, t}$ from principal component factors $p_{dim}$, which are optimization variables of the embedded maximization problem.
Worst-case demand is calculated according to 
\begin{equation*}
    c_{pca, demand, 0, t} + \sum_{dim = 1}^{n_{dim}} c_{pca, demand, dim, t} p_{dim} - P_{demand, wc, t} = 0  \ \forall t \in \mathcal{T}
\end{equation*}
$c_{pca, demand, 0, t}$ is the mean of the historical demand data, and $c_{pca, demand, dim, t}$ are the principal components that are parameters in the optimization problem and computed during preprocessing using the scikit-learn package \citep{scikit-learn}.
There are $n_{dim}$ principal components where $n_{dim}$ is the number of dimensions of the latent space.
Similarly, we compute worst-case wind and solar capacity factors according to
\begin{equation*}
    c_{pca, c, 0, t} + \sum_{dim = 1}^{n_{dim}} c_{pca, c, dim, t} p_{dim} - f_{c, wc, t} \le 0 \ \forall c \in \{solar, wind\} \ \forall t \in \mathcal{T},
\end{equation*}
where we additionally introduce the equation
\begin{equation*}
    - f_{c, wc, t} \le 0 \ \forall c \in \{solar, wind\} \ \forall t \in \mathcal{T}
\end{equation*}
to ensure that negative capacity factors due to PCA truncation errors do not occur.

Finally, we provide the equations for the convex hull.
There are two equivalent formulations of the convex hull around the historical data in the principal component space\citep{avisHowGoodAre1997}.
First, the convex hull can be encoded by an intersection of halfspaces as
\begin{equation*}
   \sum_{dim = 1}^{n_{dim}} a_{dim, j} p_{dim} + b_{conv, j} \le 0 \ \forall j \in \mathcal{J},
\end{equation*}
where $a_{dim, j}$ and $b_{conv, j}$ are parameters describing the hyperplanes bounding the historical data.
$a_{dim, j}$ and $b_{conv, j}$ are obtained using the 'ConvexHull' function of the SciPy package \citep{virtanenSciPyFundamentalAlgorithms2020}.

Second, the convex hull can be described as a linear combination of its vertices, i.e., 
\begin{equation*}
   \sum_{v = 1}^{n_{v}} \alpha_{v} c_{data, v, dim} = p_{dim} \ \forall dim \in \{1, \dots, n_{dim}\},
\end{equation*}
where $c_{data, v, dim}$ is the $dim$-th component of the latent representation vector $\bm c_{data,v}$ of the vertex $v$ and $\alpha_{v}$ are the factors describing the weights of the vertices with the restriction
\begin{equation*}
   \sum_{v = 1}^{n_{v}} \alpha_{v} = 1.
\end{equation*}
The vertices are again obtained using the 'ConvexHull' function of the SciPy package \citep{virtanenSciPyFundamentalAlgorithms2020}.
In preliminary studies, we found that representing the uncertainty set in the optimization problem as a convex combination of the vertices yielded the best computational performance in our case study.
Hence, we formulate the convex hull as a convex combination of its vertices in all our numerical experiments.

\subsection{Cost Parameters}\label{Parameters}
We adjust all price data for inflation to bring them to a $2022$ base.
To this end, we calculate an average annual producer price index \citep{eurostatDomesticProducerPrices2024} $ppi_{yr}$ from monthly data $ppi_{yr, mon}$:
\begin{align*}
    ppi_{yr} = \frac{1}{12} \sum_{mon = 1}^{12} ppi_{yr, mon}
\end{align*}
 The adjusted price $pr_{2022}$ in \si{\euro_{2022}} is then calculated from the price $pr_{yr}$ for the reference year $yr$ and the corresponding average annual producer price indices:
\begin{align*}
    pr_{2022} = pr_{yr} \frac{ppi_{2022}}{ppi_{yr}}
\end{align*}
Investment and operational costs for wind and solar are obtained from \cite{kostLevelizedCostElectricity}.
Specifically, we use the center of the provided intervals for onshore wind and utility-scale PV data.
Battery investment costs are taken as the center of the interval given by \cite{figgenerDevelopmentBatteryStorage2023} for large-scale energy storage and multiplied by the assumed energy storage ratio of \SI{4}{\kilo\watt\hour\per\kilo\watt} to obtain the battery capacity in \SI{}{\kilo\watt}, which is the unit of the battery capacity (see Section \ref{sec:Problem Formulation}).
Operational costs are assumed to be \SI{2}{\percent} of investment costs following \cite{kostLevelizedCostElectricity}.
Investment costs and variable costs of the diesel generators are calculated according to information from Spanish legislation \citep{ministeriodeindustriaenergiayturismoRealDecreto7382015}.
To calculate diesel generator cost, we assume the use of \SI{11.5}{\mega\watt} four-stroke diesel engines, as this is the size of the diesel generators most recently added to the island generation park \citep{ministeriodeindustriaenergiayturismoRealDecreto7382015}.
Variable costs are calculated according to Spanish legislation \citep{ministeriodeindustriaenergiayturismoRealDecreto7382015}:
\begin{equation*}
    c_{diesel, var} = co_{fuel} + co_{start-up} + co_{maintenance} + co_{dispatch} + co_{CO2} + co_{reduction}
\end{equation*}
Variable fuel costs $co_{fuel} = $ \SI{0.12}{\euro_{2022} \per\kilo\watt\hour} were calculated from data released by Spanish ministerial orders:
\begin{equation*}
    co_{fuel} = \frac{pr_{fuel} + pr_{logistics}}{LHV \eta_{therm}}
\end{equation*}
Here, we assumed a fuel price for \SI{1}{\percent} fuel oil in the Canaries of \SI{448.38}{\euro_{2022}\per\ton}, which is the average over the first \num{6} months of 2023 \citep{ministerioparalatransicionecologicayelretodemograficoResolucion21Septiembre2023} and logistic costs for transport to La Palma of \SI{103.47}{\euro_{2022}\per\ton} \citep{ministerioparalatransicionecologicayelretodemograficoOrdenTED13152022}.
The lower heating value for the \SI{1}{\percent} fuel oil is taken as \SI{11214.46}{\kilo\watt\hour\per\ton} \citep{ministerioparalatransicionecologicayelretodemograficoResolucion20Junio2023} and we use the average thermal efficiency of the currently used diesel generators $\eta_{therm} = 0.41$ \citep{gobiernodecanariasAnuarioEnergeticoCanarias2023}.
We neglect the dispatch start-up costs $co_{start-up}$ as we model the diesel generators in a simplified way without minimum part load constraints.
We assume the variable maintenance costs as \SI{0.042}{€_{2022}\per\kilo\watt\hour}, identical to the current \SI{11.5}{\mega\watt} four-stroke diesel engines \citep{ministeriodeindustriaenergiayturismoRealDecreto7382015}.
Dispatch regulation band costs $co_{dispatch}$ are \SI{1}{\%} of the sum of variable fuel costs and costs of dispatch emission rights.
Costs of dispatch emission rights $co_{CO2} =$ \SI{0.079}{€_{2022} \per\kilo\watt\hour} are calculated according to \citep{ministeriodeindustriaenergiayturismoRealDecreto7382015}:
\begin{equation*}
    co_{CO2} = pr_{CO2L} fa_{emission} co_{corr, p} co_{corr, e}
\end{equation*}
Here, the price of dispatch allowances $pr_{CO2L}$ is \SI{80.821}{€_{2022} \per\ton} \citep{ministerioparalatransicionecologicayelretodemograficoResolucionFebrero20232023}, the emission factor $fa_{emission}$ is \SI{0.62}{\ton_{CO2}\per\mega\watt\hour} \citep{ministeriodelapresidenciaRealDecreto13702006}, the correlation factor $co_{corr, p} $ is \num{1.028}. Pragmatically, we assume the correlation factor $co_{corr, e}$ as \num{1} since no values have been published by the responsible ministry yet.
Furthermore, the reduction of variable costs due to income or avoided costs unrelated to electricity production $co_{reduction}$ is neglected.
\subsection{Problem Formulation}\label{Problem Formulation}
In this Section, we provide the full problem formulation for the La Palma robust energy system design (see Section \si{4} in the main manuscript).
{\allowdisplaybreaks
\begin{alignat*}{3}\tag{Full Design Problem}\label{FullDesignProblem}
    &\underset{\mathbf{x} \in \mathbb{R}_{\ge 0}^{|\mathcal{C}| + 2}, \mathbf{z}_s \in \mathbb{R}^{|\mathcal{S}||\mathcal{T}|(|\mathcal{C}| + 4)}}{\min} & \quad &\sum_{c\in \mathcal{C} \cup \{battery\}} P_{c, peak} \left(\frac{c_{inv, c}}{f_{AN}} + c_{fix, c}\right) &\\&&& + 365 \sum_{s \in \mathcal{S}} \omega_s \sum_{c\in \mathcal{C}}{ \sum_{t \in \mathcal{T}} \Delta_t P_{c, peak} f_{c, s, t} c_{var, c}},&\\
    &\text{s.t.} &&P_{demand, s, t} + P_{battery, in, s, t} - P_{battery, out, s, t} - \sum_{c\in \mathcal{C}} P_{c, peak} f_{c, s, t} \le 0 \ \forall s \in \mathcal{S} \ \forall t \in \mathcal{T},&\\
    &&& - f_{c, s, t} \le 0 \ \forall c \in \mathcal{C} \ \forall s \in \mathcal{S} \ \forall t \in \mathcal{T},&\\
     &&& f_{diesel, s, t} - 1 \le 0 \ \forall s \in \mathcal{S} \ \forall t \in \mathcal{T},&\\
    &&& f_{c, s, t} -  c_{cap, c, s, t} \le 0 \ \forall c \in \{wind, solar\} \ \forall s \in \mathcal{S} \ \forall t \in \mathcal{T},&\\
    &&& P_{battery, out, s, t} - P_{battery, peak} \le 0 \ \forall s \in \mathcal{S} \ \forall t \in \mathcal{T},&\\
    &&& - P_{battery, in, s, t} \le 0 \ \forall s \in \mathcal{S} \ \forall t \in \mathcal{T},&\\
    &&& P_{battery, in, s, t} - P_{battery, peak} \le 0 \ \forall s \in \mathcal{S} \ \forall t \in \mathcal{T},&\\
    &&& - P_{battery, out, s, t} \le 0 \ \forall s \in \mathcal{S} \ \forall t \in \mathcal{T},&\\
    &&& E_{battery, s, t} - E_{battery, peak} \le 0 \ \forall s \in \mathcal{S} \ \forall t \in \mathcal{T},&\\
    &&& - E_{battery, s, t} \le 0 \ \forall s \in \mathcal{S} \ \forall t \in \mathcal{T},&\\
    &&& \dot E_{battery, s, t} - \frac{1}{3600} (\eta_{in} P_{battery, in, s, t} - \frac{1}{\eta_{out}} P_{battery, out, s, t}) = 0 \ \forall s \in \mathcal{S}
    \ \forall t \in \mathcal{T},&\\
    &&& \dot E_{battery, s, t} \Delta t - (E_{battery, s, t} - E_{battery, s, t-1}) = 0 \ \forall s \in \mathcal{S} \ \forall t \in \mathcal{T},&\\
    &&& E_{battery, s, 0} - 0.5 E_{battery, peak} = 0 \ \forall s \in \mathcal{S},&\\
    &&& E_{battery, peak} - 4 P_{battery, peak}  = 0, &\\
    &&&\underset{\mathbf{y}\in \mathcal{Y}}{\max}\underset{\mathbf{z}\in \mathcal{Z}_{epi}(\mathbf{x}, \mathbf{y})}{\min} \quad e_{epi} \le 0,&\\
\end{alignat*}
with 
\begin{alignat*}{1}
    &\mathbf{x} =  [P_{1, peak}, \dots, P_{|\mathcal{C}|, peak}, P_{battery, peak}, E_{battery, peak}]\\
    &\mathbf{z}_s =  [P_{battery, in, 1, 1}, \dots, P_{battery, in, 1, |\mathcal{T}|}, P_{battery, in, 2, 1}, \dots, P_{battery, in, |\mathcal{S}|, |\mathcal{T}|}, P_{battery, out, 1, 1}, \dots, P_{battery, out, 1, |\mathcal{T}|},\\
    &P_{battery, out, 2, 1}, \dots, P_{battery, out, |\mathcal{S}|, |\mathcal{T}|}, E_{battery, 1, 1}, \dots, E_{battery, 1, |\mathcal{T}|}, E_{battery, 2, 1}, \dots, E_{battery, |\mathcal{S}|, |\mathcal{T}|}, \dot E_{battery, 1, 1}, \dots,\\
    &\dot E_{battery, 1, |\mathcal{T}|}, \dot E_{battery, 2, 1}, \dots, \dot E_{battery, |\mathcal{S}|, |\mathcal{T}|}, f_{1, 1, 1}, \dots, f_{1, 1, |\mathcal{T}|}, f_{1, 2, 1}, \dots, f_{1, |\mathcal{S}|, |\mathcal{T}|}, f_{2, 1, 1}, \dots, f_{|\mathcal{C}|, |\mathcal{S}|, |\mathcal{T}|}]\\
    &\mathbf{y} =  [f_{wind, wc, 1}, \dots, f_{wind, wc, |\mathcal{T}|}, f_{solar, wc, 1}, \dots, f_{solar, wc, |\mathcal{T}|}, p_{1}, \dots, p_{|n_{dim}|}, P_{demand, wc, 1}, \dots, P_{demand, wc, |\mathcal{T}|},\\
    &\alpha_{1}, \dots, \alpha_{|n_{v}|}]\\
    &\mathcal{Y} = \{\mathbf{y} \in \mathbb{R}^{3|\mathcal{T}|+ |n_{v}| + |n_{dim}|} \ | \ \mathbf{g_{y}}(\mathbf{y}) \le \mathbf{0} \wedge \mathbf{h_{y}}(\mathbf{y}) = \mathbf{0}\},\\
    &\mathbf{g_{y}}(\mathbf{y}):\\
    &c_{pca, c, 0, t} + \sum_{dim = 1}^{n_{dim}} c_{pca, c, dim, t} p_{dim} - f_{c, wc, t} \le 0 \ \forall c \in \{solar, wind\} \ \forall t \in \mathcal{T},\\
    & - f_{c, wc, t} \le 0 \ \forall c \in \{solar, wind\} \ \forall t \in \mathcal{T}.\\
    &\mathbf{h_{y}}(\mathbf{y}):\\
    & c_{pca, c, 0, t} + \sum_{dim = 1}^{n_{dim}} c_{pca, c, dim, t} p_{dim} - P_{demand, wc, t} = 0 \ \forall c \in \{demand\} \ \forall t \in \mathcal{T},\\
    & \sum_{v = 1}^{n_{v}} \alpha_{v} c_{data, v, dim} = p_{dim} \ \forall dim \in \{1, \dots, n_{dim}\},\\
    & \sum_{v = 1}^{n_{v}} \alpha_{v} = 1.\\
    &\mathbf{z} = [P_{battery, in, wc, 1}, \dots, P_{battery, in, wc, |\mathcal{T}|}, P_{battery, out, wc, 1}, \dots, P_{battery, out, wc, |\mathcal{T}|}, E_{battery, wc, 1}, \dots, E_{battery, wc, |\mathcal{T}|},\\
    &\dot E_{battery, wc, 1}, \dots, \dot E_{battery, wc, |\mathcal{T}|}, f_{diesel, wc, 1}, \dots, f_{diesel, wc, |\mathcal{T}|, e_{epi}}]\\
    &\mathcal{Z}_{epi}(\mathbf{x}, \mathbf{y}) = \{ \mathbf{z} \in \mathbb{R}^{5|\mathcal{T}| + 1} \ | \ \mathbf{g}_{\mathbf{z}}(\mathbf{x}, \mathbf{z}) \le \mathbf{0} \wedge \mathbf{h}_{\mathbf{z}}(\mathbf{x}, \mathbf{z}) = \mathbf{0} \wedge \\
    & P_{demand, wc, t} + P_{battery, in, wc, t} - P_{battery, out, wc, t} - \sum_{c\in \mathcal{C}} P_{c, peak} f_{c, wc, t} - e_{epi} \le 0 \ \forall t \in \mathcal{T}\}\\
    &\mathbf{g}_{\mathbf{z}}(\mathbf{x}, \mathbf{z}):\\
    & P_{battery, out, s, t} - P_{battery, peak} \le 0 \ \forall t \in \mathcal{T},\\
    & - P_{battery, in, s, t} \le 0 \ \forall t \in \mathcal{T},\\
    & P_{battery, in, s, t} - P_{battery, peak} \le 0 \ \forall t \in \mathcal{T},\\
    & - P_{battery, out, s, t} \le 0 \ \forall t \in \mathcal{T},\\
    & E_{battery, s, t} - E_{battery, peak} \le 0 \ \forall t \in \mathcal{T},\\
    & - E_{battery, s, t} \le 0\ \forall t \in \mathcal{T}.\\
    &\mathbf{h}_{\mathbf{z}}(\mathbf{x}, \mathbf{z}): \\
    &\dot E_{battery, wc, t} - \frac{1}{3600} (\eta_{in} P_{battery, in, wc, t} - \frac{1}{\eta_{out}} P_{battery, out, wc, s, t}) = 0 \ \forall t \in \mathcal{T},\\
    & E_{battery, s, 0} - 0.5 E_{battery, peak} = 0,\\   
    & \dot E_{battery, s, t} \Delta t - (E_{battery, s, t} - E_{battery, s, t-1}) = 0 \ \forall t \in \mathcal{T},\\
    & f_{diesel, wc, t} - 1 = 0 \ \forall t \in \mathcal{T}.\\
\end{alignat*}
}

Explanations for all symbols were introduced in Sections \ref{sec:Problem Formulation} and \ref{Parameters}.
Note the presence of a coupling equality constraint, i.e,  $E_{battery, s, 0} - 0.5 E_{battery, peak} = 0$, which we did not remove since no convergence issues occurred.
Furthermore, note the independence of the energy system model equations $\mathbf{g}_{\mathbf{z}}(\mathbf{x}, \mathbf{z})$ and $\mathbf{h}_{\mathbf{z}}(\mathbf{x}, \mathbf{z})$ of $\mathbf{y}$.

\subsection{Lifted Problem Formulation}
We apply the lifting approach described in Section \si{3} of the main manuscript to transform the MAXMIN problem into a single-level NLP.
{\allowdisplaybreaks
\begin{alignat*}{4}
    & \underset{\mathbf{y}, \mathbf{z}, \boldsymbol{\lambda}, \boldsymbol{\mu}}{\max}  &\qquad  & \mathcal{L}_{ll}(\mathbf{x},\mathbf{y},\mathbf{z}, \boldsymbol{\lambda}, \boldsymbol{\mu}),&\\
    & \text{s.t.} &&\nabla_{\mathbf{z}} \mathcal{L}_{ll}(\mathbf{x},\mathbf{y},\mathbf{z}, \boldsymbol{\lambda}, \boldsymbol{\mu}) = \mathbf{0},&\\
    &&& \boldsymbol{\mu}\ge \mathbf{0},&\\
    &&& \boldsymbol{\mu}^\intercal\mathbf{g}_{\mathbf{z}}(\mathbf{x},\mathbf{y},\mathbf{z}) = \mathbf{0},&\\
    &&&\mathbf{g}_{\mathbf{z}}(\mathbf{x}, \mathbf{y}, \mathbf{z}) \le \mathbf{0},&\\
    &&&\mathbf{h}_{\mathbf{z}}(\mathbf{x}, \mathbf{z}) = \mathbf{0},&\\
    &&& c_{pca, c, 0, t} + \sum_{dim = 1}^{n_{dim}} c_{pca, c, dim, t} p_{dim} - f_{c, wc, t} \le 0 \ \forall c \in \{solar, wind\} \ \forall t \in \mathcal{T},&\\
    &&&- f_{c, wc, t} \le 0 \ \forall c \in \{solar, wind\} \ \forall t \in \mathcal{T},&\\
     &&&c_{pca, c, 0, t} + \sum_{dim = 1}^{n_{dim}} c_{pca, c, dim, t} p_{dim} - P_{demand, wc, t} = 0 \ \forall c \in \{demand\} \ \forall t \in \mathcal{T},&\\
    &&& \sum_{v = 1}^{n_{v}} \alpha_{v} c_{data, v, dim} = p_{dim} \ \forall dim \in \{1, \dots, n_{dim}\},&\\
    &&& \sum_{v = 1}^{n_{v}} \alpha_{v} = 1 ,&\\
\end{alignat*}
with 
\begin{alignat*}{1}
    &\mathbf{y} =  [f_{wind, wc, 1}, \dots, f_{wind, wc, |\mathcal{T}|}, f_{solar, wc, 1}, \dots, f_{solar, wc, |\mathcal{T}|}, p_{1}, \dots, p_{|n_{dim}|}, P_{demand, wc, 1}, \dots, P_{demand, wc, |\mathcal{T}|},\\
    &\alpha_{1}, \dots, \alpha_{|n_{v}|}] \in \mathbb{R}^{3|\mathcal{T}|+ |n_{v}| + |n_{dim}|},\\
    &\mathbf{z} = [P_{battery, in, wc, 1}, \dots, P_{battery, in, wc, |\mathcal{T}|}, P_{battery, out, wc, 1}, \dots, P_{battery, out, wc, |\mathcal{T}|}, E_{battery, wc, 1}, \dots, E_{battery, wc, |\mathcal{T}|},\\
    &\dot E_{battery, wc, 1}, \dots, \dot E_{battery, wc, |\mathcal{T}|}, f_{diesel, wc, 1}, \dots, f_{diesel, wc, |\mathcal{T}|, e_{epi}}] \in \mathbb{R}^{5|\mathcal{T}| + 1} \\
    &\boldsymbol{\lambda} \in \mathbb{R}^{3|\mathcal{T}| + 1},\\
    &\boldsymbol{\mu} \in \mathbb{R}^{7|\mathcal{T}|},\\
    &\mathbf{g}_{\mathbf{z}}(\mathbf{x}, \mathbf{y}, \mathbf{z}): \\
    &P_{demand, wc, t} + P_{battery, in, wc, t} - P_{battery, out, wc, t} - \sum_{c\in \mathcal{C}} P_{c, peak} f_{c, wc, t} - e_{epi} \le 0 \ \forall t \in \mathcal{T},\\
    & P_{battery, out, wc, s, t} - P_{battery, peak} \le 0 \ \forall t \in \mathcal{T},\\
    & - P_{battery, in, wc, s, t} \le 0 \ \forall t \in \mathcal{T},\\
    & P_{battery, in, wc, s, t} - P_{battery, peak} \le 0 \ \forall t \in \mathcal{T},\\
    & - P_{battery, out, wc, s, t} \le 0 \ \forall t \in \mathcal{T},\\
    & E_{battery, wc, s, t} - E_{battery, peak} \le 0 \ \forall t \in \mathcal{T},\\
    & - E_{battery, wc, s, t} \le 0\ \forall t \in \mathcal{T}.\\
    &\mathbf{h}_{\mathbf{z}}(\mathbf{x}, \mathbf{z}): \\
    &\dot E_{battery, wc, t} - \frac{1}{3600} (\eta_{in} P_{battery, in, wc, t} - \frac{1}{\eta_{out}} P_{battery, out, wc, s, t}) = 0 \ \forall t \in \mathcal{T},\\
    & E_{battery, wc, s,0} - 0.5 E_{battery, peak} = 0,\\
    & \dot E_{battery, wc, s, t} \Delta t - (E_{battery, wc, s, t} - E_{battery, wc, s, t-1}) = 0 \ \forall t \in \mathcal{T},\\
    & f_{diesel, wc, t} - 1 = 0 \ \forall t \in \mathcal{T}.\\
    &\mathcal{L}_{ll}(\mathbf{x},\mathbf{y},\mathbf{z}, \boldsymbol{\lambda}, \boldsymbol{\mu}) = e_{epi} + \boldsymbol{\lambda}^\intercal \mathbf{h}_{\mathbf{z}}(\mathbf{x}, \mathbf{z}) + \boldsymbol{\mu}^\intercal \mathbf{g}_{\mathbf{z}}(\mathbf{x},\mathbf{y},\mathbf{z})
\end{alignat*}
}
Where $\mathcal{L}_{ll}$ is the Lagrangian function and $\boldsymbol{\mu}$ and $\boldsymbol{\lambda}$ are the Lagrange multipliers for the inequalities and equalities, respectively.
Explanations for all remaining symbols were introduced in Sections \ref{sec:Problem Formulation} and \ref{Parameters}.

\section{Solver Settings}\label{Optimization Parameters}
We describe the libDIPS \citep{zinglerLibDIPSDiscretizationBasedSemiInfinite2023} settings as well as the Gurobi (version 11.0.3) \citep{gurobi} settings that we used to generate the results shown in Section 3 and 4 of the main manuscript.
libDIPS was used for the solution of hierarchical programs and Gurobi was used to solve the problems formulated by libDIPS.
\begin{table}[h]
    \begin{center}
        \begin{tabular}{c|c}
            \hline
            Name & Value\\
            \hline
            \multicolumn{2}{c}{sip\_bnf\_demo.cpp} \\
            \hline
            feas\_tol & \num{5e-2} \\ 
            abs\_tol\_lbp & \num{5e-3} \\
            rel\_tol\_lbp & \num{5e-3} \\
            abs\_tol\_lbp\_llp & \num{1e-2} \\
            rel\_tol\_lbp\_llp & \num{5e-3} \\
            max\_time & \num{10800} \\
            \hline
            \multicolumn{2}{c}{esip\_bnf\_demo.cpp} \\
            \hline
            feas\_tol & \num{5e-2} \\
            \hline
            \multicolumn{2}{c}{esip\_bnf\_solver.hpp} \\
            \hline
            abs\_tol & \num{1e-2} \\
            rel\_tol & \num{2e-1} \\
            abs\_tol\_minmax & \num{2.5e-2} \\
            rel\_tol\_minmax & \num{1e-3} \\
            \hline
            \multicolumn{2}{c}{Gurobi\_solver.hpp} \\
            \hline
            non\_convex & \num{2} \\
            num\_threads & \num{4} \\
            OBBT & \num{3} \\
            Presolve & \num{2} \\
            FuncNonlinear & \num{1} \\
            \hline
        \end{tabular}
        \caption{libDIPS \citep{zinglerLibDIPSDiscretizationBasedSemiInfinite2023} settings.\label{tab:dipssettings}}
    \end{center}
\end{table}
Table \ref{tab:dipssettings} shows the settings that deviate from the default  settings of libDIPS and Gurobi. If not otherwise stated, the default values were used. 
The filenames refer to the files where these settings were changed, e.g., 'sip\_bnf\_demo.cpp'.

The libDIPS parameters that were changed to solve the SIP problem are:
'feas\_tol' is the feasibility tolerance of the semi-infinite constraint for which the problem is considered feasible.
'abs\_tol\_lbp', 'abs\_tol\_llp', 'rel\_tol\_lbp', and 'rel\_tol\_llp' are the absolute and relative optimality tolerances for the lower bounding and lower-level problems, respectively.
'max\_time' is the maximum wall-clock time before the optimization is aborted if convergence is not reached.

The libDIPS parameters that were changed to solve the ESIP problem are:
'feas\_tol' is the feasibility tolerance of the semi-infinite existence constraint for which the problem is considered feasible.
'abs\_tol', and 'rel\_tol' are the absolute and relative optimality tolerances for the lower bounding problem.
'abs\_tol\_minmax', and 'rel\_tol\_minmax' are the absolute and relative optimality tolerances for the min-max implementation in libDIPS and are multiplied with respective factors in the code to provide tolerances for the medial- and lower-level problems.

Finally, the altered parameters specific to the Gurobi \citep{gurobi} solver are given:
'non\_convex' enables Gurobi to handle nonconvex problems.
'num\_threads' informs Gurobi about the number of threads to use in parallel to solve the optimization problem.
'OBBT' sets the amount of work allowed for optimality-based bound tightening, $3$ is the most agressive setting.
'Presolve' sets the presolve level, $2$ is agressive presolving, taking more time but yielding a tighter model.
'FuncNonlinear' determines how nonlinear equations are handled. $1$ approximates nonlinear equations by outer approximation.

\section*{Nomenclature}
This Section extends the Nomenclature of the main manuscript.
Throughout the manuscript, scalar-valued quantities are denoted in regular font, e.g., $x$, vector-valued quantities are denoted in bold font, e.g., $\mathbf{x}$, and set-valued quantities are denoted in calligraphic font, e.g., $\mathcal{X}$.
\subsection*{Abbreviations}
\begin{table}[H]
\begin{tabular}{ll}
    ESIP & existence-constrained semi-infinite program \\
    MILP & mixed-integer linear program \\
    MW & megawatt \\
    NLP & nonlinear program \\
    LHV & lower heating value \\
    RESD & robust energy system design \\
    PC & principal component \\
    PCA & principal component analysis \\
    PV & photo voltaic \\
    SIP & semi-infinite program \\
    TAC & total annualized cost \\
\end{tabular}
\end{table}

\subsection*{Greek Symbols}
\begin{table}[H]
\begin{tabular}{ll}
    $\alpha$ & weighting factor \\
    $\Delta_t$ & time-step length \\
    $\eta$ & efficiency factor \\
    $\boldsymbol{\lambda}$ & Lagrange multipliers for equality constraints \\
    $\boldsymbol{\mu}$ & Lagrange multipliers for inequality constraints \\
    $\omega_s$ & weighting factor for representative scenarios \\
\end{tabular}
\end{table}

\subsection*{Latin Symbols}
{\allowdisplaybreaks
\begin{table}[H]
\begin{tabular}{ll}
    $a$ & factor in convex hull equation \\
    $b_{ll}$ & lower-level binary variable for on/off decision \\
    $b_{conv}$ & factor in convex hull equation \\
    $b_u$ & upper-level binary variable for on/off decision \\
    $b_m$ & medial-level binary variable to reformulate the maximum in the objective function\\
    $c_{cap}$ & capacity factor parameter \\
    $c_{fix}$ & fixed costs \\
    $c_{inv}$ & investment costs \\
    $c_{pca}$ & principal component analysis factor \\
    $c_{var}$ & variable costs \\
    $co_{corr, e}$ & correlation factor \\
    $co_{corr, p}$ & correlation factor \\
    $co_{CO_2}$ & cost of dispatch emission rights\\
    $co_{dispatch}$ & dispatch regulation band costs\\
    $co_{fuel}$ & variable fuel costs\\
    $co_{maintenance}$ & variable maintenance costs\\
    $co_{start-up}$ & dispatch start-up costs\\
    $co_{reduction}$ & cost coefficient for reduced costs\\
    $dim$ & dimensionality of the latent space \\
    $e$ & energy supply gap \\
    $E$ & energy \\
    $\dot{E}$ & rate of charge \\
    $e_{epi}$ & auxiliary variable for energy supply gap \\
    $f$ & capacity factor \\
    $fa_{emission}$ & emission factor \\
    $f_{AN}$ & annuity factor \\
    $\mathbf{g}$ & vector of inequality constraints \\
    $h_{hub}$ & wind turbine hub height \\
    $\mathbf{h}$ & vector of equality constraints \\
    \end{tabular}
    \end{table}
    \begin{table}[H]
    \begin{tabular}{ll}
    $i$ & interest \\
    $I$ & irradiance \\
    $\mathcal{J}$ & index set \\
    $\mathcal{K}$ & discretization set \\
    $llpobj$ & auxiliary variable for lower-level problem objective \\
    $\mathcal{L}$ & Lagrangian function \\
    $mlpobj$ & auxiliary variable for medial-level problem objective \\
    $\mathcal{N}$ & discretization set \\
    $n$ & number of a quantity \\
    $p$ & principal component \\
    $P$ & power \\
    $ppi$ & producer price index \\
    $pr$ & price \\
    $pr_{CO2L}$ & price for carbon dioxide emissions \\
    $pr_{fuel}$ & fuel price \\
    $pr_{logistics}$ & logistic costs \\
    $\mathcal{S}$ & set of representative scenarios \\
    $T$ & length of time horizon \\
    $t$ & time \\
    $\mathcal{T}$ & set of time steps \\
    $v_{wind}$ & wind speed\\
    $x_1$ & design variable for component 1\\
    $x_2$ & design variable for component 2\\
    $\mathbf{x}$ & vector of design variables \\
    $\mathcal{X}$ & feasible set of design variables \\
    $y$ & uncertain demand\\
    $\mathbf{y}$ & vector of uncertain variables \\
    $\mathcal{Y}$ & feasible set of uncertain variables\\
    $z$ & operational variable \\
    $z_0$ & roughness length\\
    \end{tabular}
    \end{table}
    \begin{table}[H]
    \begin{tabular}{ll}
    $z_{1, u}$ & upper-level operational variable for component 1\\
    $z_{1, ll}$ & lower-level operational variable for component 1\\
    $\mathbf{z}$ & vector of operational variables \\
    $\mathcal{Z}$ & feasible set of operational variables \\
\end{tabular}
\end{table}

\subsection*{Subscripts}
\begin{table}[H]
\begin{tabular}{ll}
    $battery$ & battery \\
    $c$ & component \\
    $data$ & data points \\
    $days$ & days \\
    $demand$ & demand \\
    $diesel$ & diesel \\
    $dim$ & dimension \\
    $disc$ & discretization \\
    $e$ & energy supply gap \\
    $en$ & energy system \\
    $epi$ & epigraph reformulation \\
    $eq$ & equations \\
    $fix$ & fixed \\
    $hub$ & hub height \\
    $in$ & input \\
    $j$ & index \\
    $k$ & index \\
    $l$ & index \\
    $mon$ & month \\
    $nom$ & nominal \\
    $out$ & output \\
    $p$ & principal components \\
    \end{tabular}
    \end{table}
    \begin{table}[H]
    \begin{tabular}{ll}
    $pca$ & principal component analysis \\
    $peak$ & peak power \\
    $ref$ & reformulation \\
    $s$ & representative scenarios \\
    $solar$ & solar PV \\
    $t$ & time step \\
    $therm$ & thermal \\
    $v$ & vertices \\
    $y$ & upper level \\
    $var$ & variable \\
    $wc$ & worst-case \\
    $wind$ & wind turbine \\
    $x$ & design variables \\
    $y$ & uncertain variables \\
    $yr$ & year \\
    $\mathbf{z}$ & operational variables \\
    $10m$ & \si{10} meters \\
\end{tabular}
\end{table}

\bibliographystyle{apalike}
  \renewcommand{\refname}{References}  
  \bibliography{bibliography.bib}